\documentclass[journal]{IEEEtran}
\usepackage{tikz}
\usetikzlibrary{calc}
\usepackage[compress]{cite}
\usetikzlibrary{automata,positioning}
\usetikzlibrary{arrows}
\definecolor {processblue}{cmyk}{0,0,0,0.17}
\definecolor{light-gray}{gray}{0.85}
\usepackage{minibox}
\usepackage{caption}
\usepackage{relsize}
\usepackage{fixltx2e}
\usepackage{url}
\usepackage{multirow}
\usepackage{float}
\usepackage{graphicx}
\usepackage[cmex10]{amsmath}
\usepackage[font={small}]{caption}
\usepackage{amssymb}
\usepackage{mathtools}
\usepackage{xfrac}

\usepackage{epstopdf}

\usepackage{xfrac}
\usepackage{multirow}
\usepackage{fixltx2e}
\usepackage{url}
\usepackage{multirow}
\usepackage{booktabs}
\usepackage{rotating}
\usepackage{array}
\usepackage{bigstrut}
\usepackage[labelformat=simple]{subcaption}

\usepackage{graphicx,times}
\usepackage{algorithm}
\usepackage{algorithmic}
\usepackage{wrapfig}
\usepackage{url}
\usepackage{amssymb}
\usepackage{paralist}
\usepackage{graphicx}
\usepackage[justification=centering]{caption}
\usepackage{tabularx}
\usepackage{subcaption}
\usepackage{xparse}
\usepackage{soul}

\newcommand*{\grayemph}[1]{%
  \tikz[baseline=(X.base)] \node[rectangle, fill=light-gray, rounded corners, inner sep=1.8mm] (X) {#1};%
}

\NewDocumentCommand{\ceil}{s O{} m}{%
  \IfBooleanTF{#1} 
    {\left\lceil#3\right\rceil} 
    {#2\lceil#3#2\rceil} 
}
\makeatletter
\newcommand*{\rom}[1]{\expandafter\@slowromancap\romannumeral #1@}
\makeatother
\ifCLASSINFOpdf
\else
\fi
\begin{document}
%
\title{Unsplittable Load Balancing in a Network of Charging Stations Under QoS Guarantees}
\author{I.~Safak~Bayram,~\IEEEmembership{Member,~IEEE,} George~Michailidis,~\IEEEmembership{Member,~IEEE,}
and~ Michael~Devetsikiotis,~\IEEEmembership{Fellow,~IEEE}
\vspace{-20 pt}
\thanks{Manuscript received 29 March 2014; revised 13 April 2014 and 30 July 2014; accepted 6 September 2014. Paper no TSG-00285-2014}
\thanks{I.~Safak Bayram (corresponding author) is with the Department
of Electrical and Computer Engineering, Texas A\&M University at Qatar, Doha,
Qatar.}
\thanks{George Michailidis
is with the Departments of Statistics and EECS, University of Michigan, Ann Arbor, MI, 48109-1107, USA.}
\thanks{ Michael Devetsikiotis is with the Department
of Electrical and Computer Engineering, North Carolina State University, Raleigh,
NC, 27695-7911 USA.}
\thanks{Emails: islam.bayram@qatar.tamu.edu, mdevets@ncsu.edu, gmichail@umich.edu}}

\maketitle
\begin{abstract}

The operation of the power grid is becoming more stressed, due to the addition of new large loads represented by Electric
Vehicles (EVs) and a more intermittent supply due to the incorporation of renewable sources. As a consequence,
the coordination and control of projected EV demand in a network of fast charging stations becomes a critical and challenging problem.
 In this paper, we introduce a game theoretic based decentralized control mechanism
to alleviate negative impacts from the EV demand. The proposed mechanism takes into consideration the
non-uniform spatial distribution of EVs that  induces uneven power demand at each charging facility, and aims to:
(i) avoid straining grid resources by offering price incentives so that customers accept being routed to less busy stations,
(ii) maximize total revenue by serving more customers with the same amount of grid resources, and
(iii) provide charging service to customers with a certain level of Quality-of-Service (QoS), the latter defined as the long term
customer blocking probability. We examine three scenarios of increased complexity that gradually approximate real world settings.
The obtained results show that the proposed framework leads to substantial performance improvements in terms of the aforementioned goals,
when compared to current state of affairs.

\end{abstract}
\noindent
\vspace{-5pt}
\section{Introduction}
\label{intro} Transportation electrification offers solutions to
an array of current societal issues, ranging from unstable oil
prices to environmental concerns. It is anticipated that a
significant proportion of vehicles (at least $10$\% of U.S. national
fleet by $2020$~\cite{jsac} and similar targets are set in Europe~\cite{revision2}) will be electricity powered
over the next decade. On the other hand, the all-electric-range
of current EVs is relatively short compared to gas powered
competitors, while the need for longer travel ranges requires
development of a network of public fast charging facilities.
However, the power grid is not well designed to accommodate the
projected EV charging demand and straining the grid beyond its
limits could easily lead to cascading failures and outages. For
instance, during peak hours the concurrent charging of only two
level-$2$ chargers can cause a failure of distribution
transformers in residential neighborhoods~\cite{jsac}.
Considering the typical amount of energy transfer to fill up an
EV battery (about $24$kWh in $30$ minutes and equivalent to the
daily average consumption of two households), a key issue becomes how to
design their architecture so as to strengthen the EV adoption
rate. Obviously, part of a long term solution that accommodates a large share
of EVs involves expansion of charging station infrastructure, but nevertheless
the design issues addressed in this study would still be critical.
 In this paper we propose a decentralized load balancing
framework in a network of fast charging stations. Apart from current literature, due to short charging sessions the proposed framework assumes that EV demand is unsplittable and optimally assigns EVs to hard-capacitated charging stations. To the best of our knowledge, little attention has been paid to investigate this problem. The overview of the proposed
system is presented in Fig.~\ref{overview2}.

Note that a network of charging facilities needs to be spatially distributed~\cite{jsac,sgc12} both due to operational limitations (e.g., transformer rating, line capacity), as well as to meet a spatially distributed EV demand. However, the two may not be perfectly aligned, thus introducing the need that
the network operator incentivizes customers to be rerouted to neighboring, but less congested stations. The most effective incentive comes in the form of
lower electricity prices charged to drivers. In our framework, we map the behavior of EVs and the system operator into a game theoretic model. The objective of EV drivers is to maximize their utility (level of satisfaction) that depends on the availability of charging service, the price paid to get the service, as well as the duration needed to complete the service request. On the other hand, the station operators want to serve as many customers as possible with the same amount of grid resources. Next, we discuss the key contributions of this paper:

\begin{figure}[t]
 \centering
 \includegraphics[width=0.95\columnwidth]{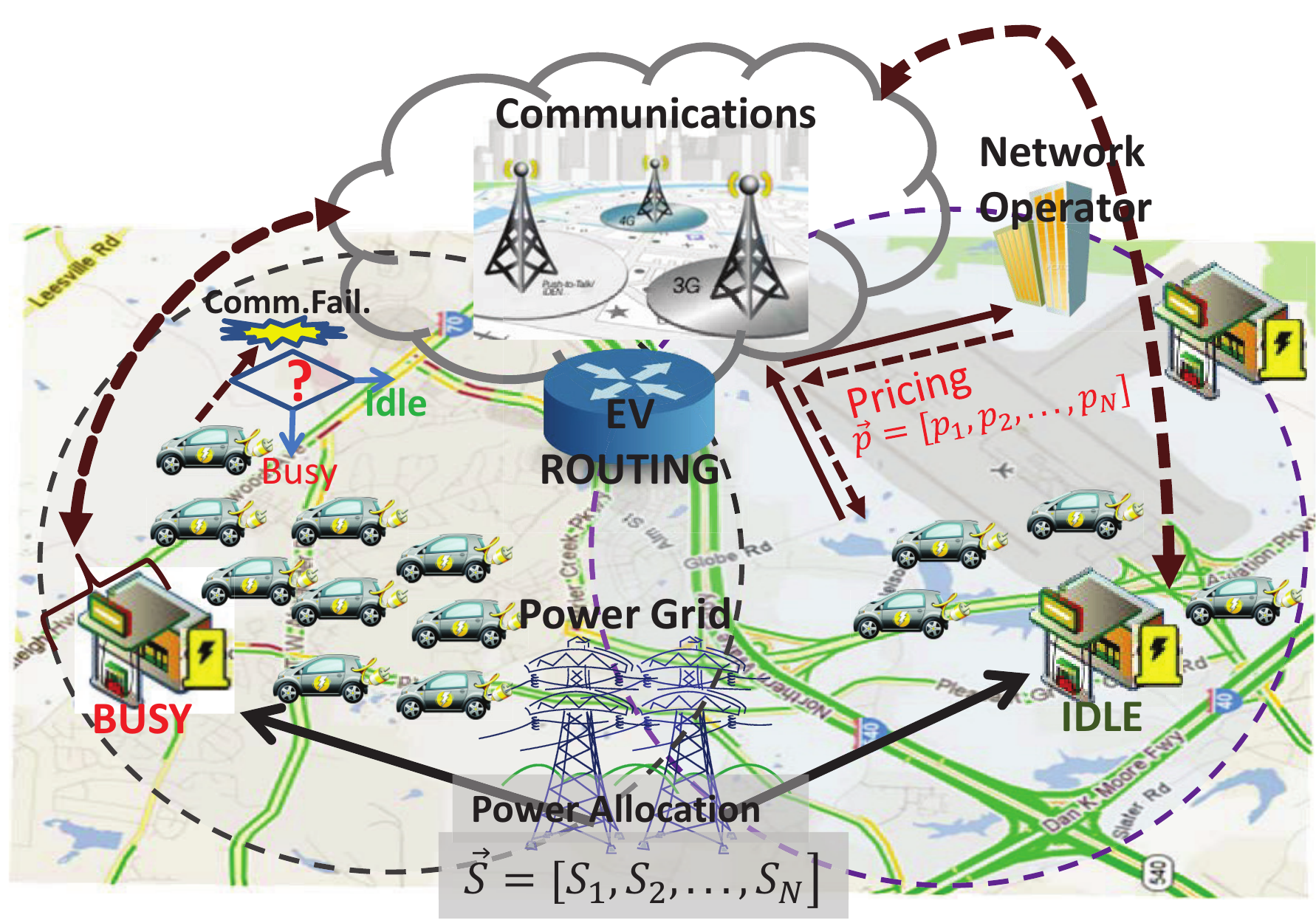}
  \caption{Overview of the paper.}\label{overview2}
 \vspace{-15pt}
\end{figure}
\begin{itemize}
  \item In a network of fast charging stations we propose a power resource allocation framework that meets QoS targets and considers the spatio-temporal EV distribution. We also introduce an Electric Vehicle Admission Control mechanism at each station that employs \emph{congestion pricing} to shape excessive demand.
  \item We introduce utility functions to capture the behavior of individual EVs and the network operator for varying network conditions.
  \item We introduce a game theoretic framework in which customer routing and load balancing are achieved through a pricing mechanism. The latter's strategy corresponds to maximizing its profit, while the strategy of the EVs is to minimizing the associated charging cost.

\end{itemize}

\section{Related Work}
\begin{figure}[t]
\centering
\begin{tikzpicture}[font=\footnotesize,edge from parent/.style={draw,thick},text width=1.5cm]
\tikzstyle{level 1}=[level distance=15mm,sibling distance=35mm]
\tikzstyle{level 2}=[level distance=15mm,sibling distance=20mm]
\node(0){Related Work}
child{node{Control}
child{node{\grayemph{Decentralized Control
\cite{HiskensX,decentralized5,revision1,price1,price3}}}}
child{node{Central Control \cite{scaglione,jsac,CallawayX,langTong1,revision3}}}}
child{node {Station Design}
child{node{Power Eng.\cite{10,canada}}}
child{node {\grayemph{Queuing Model\cite{EVaustin,jsac,queue2}}}}
};
\end{tikzpicture}
\caption{The key components of this work  correspond to the  highlighted topics.}\label{relatedWork}
\vspace{-15 pt}
\end{figure}
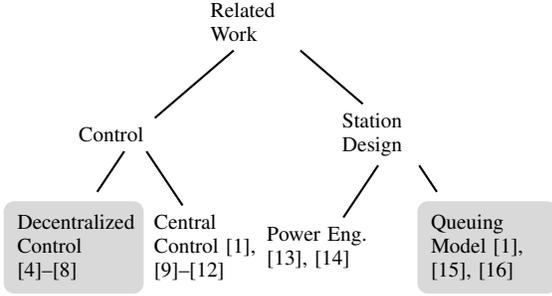
\subsection{Comparative Literature Analysis}
There has been an increasing body of literature on control
mechanisms for EV charging,
while the literature on designs for charging stations is rather sparse.
The predominant assumption for the control of EV charging is
that the vehicles are located either in their drivers' residencies or at large parking
lots and also depending on the charger technology the duration of the charging session is of the order of hours.
The vast majority of the literature threats EVs as ``smart loads" and management of their demand is achieved by
optimally adjusting the charging current. Overall, related literature can be classified into two categories: centralized and decentralized (distributed)
 control \cite{revision2}. In centralized control a network operator (dispatcher) possesses global information about all users and
to a large extent controls and mandates charging times, rates etc. For instance, the authors of \cite{scaglione} propose a direct load
management for EV chargings. Similarly, in our previous work~\cite{jsac}, we introduce a central EV allocation framework
for a network of fast charging stations.

In decentralized control, individuals choose their own service patterns based on cost minimization principles. In this case,
utilities aim to shape demand profiles by incentivizing customers via pricing-based control mechanisms. The work presented
in \cite{HiskensX} proposes a game theoretic framework to optimally control the demand of a large scale stationary EV population. Further,
the authors of \cite {decentralized5} and~\cite{revision1} present decentralized control mechanisms for EVs and a detailed survey can be found in \cite{comparativeAnalysis}. Note
that decentralized control eliminates the need for advanced
monitoring tools, whereas centralized control leads to better
utilization of charging resources.

Our proposed framework has the following distinctive features. We assume that the EVs are mobile and the charging session at fast stations takes around $25$-$30$ minutes. Hence, the primary goal of the EVs is to receive service as quickly as possible and extend their driving range. Moreover, the demand of each individual EV is assumed to be unsplittable and should be served at exactly one station. Also the load of each station should be kept below its serving capacity (power drawn from the grid plus the local energy storage unit) to avoid disruptions. To that end, the nature of the problem requires a decentralized load balancing solution for unsplittable EV demand at hard capacitated charging facilities. Obviously the routing capabilities depend on the load at each station. If the customer demand increases substantially, station operators need to upgrade their storage capabilities and/or increase the power drawn from the grid. However, the proposed framework enables station operators to defer power grid upgrades, which would occur over a broader time horizon
and  in a more cost effective way.

Related work on charging station design can be classified into two categories.
The first category includes queueing based models, where the goal is to evaluate the
charging station performances with respect to long-term statistical metrics e.g., blocking probability, waiting time, etc.
For example, the work in \cite{EVaustin} proposes an M/M/s queueing
based mathematical model for EV demand for fast charging
stations located near highway exits. On the other hand, the design perspective of the second category is related to
power engineering. The authors of \cite{10} propose a charging station architecture
with a DC bus system, and similar to our model they employ a
local energy storage unit to alleviate the stress in the power
grid. Related literature is presented in Fig. \ref{relatedWork}.

\begin{figure}[t]
 \centering
 \includegraphics[width=0.75\columnwidth]{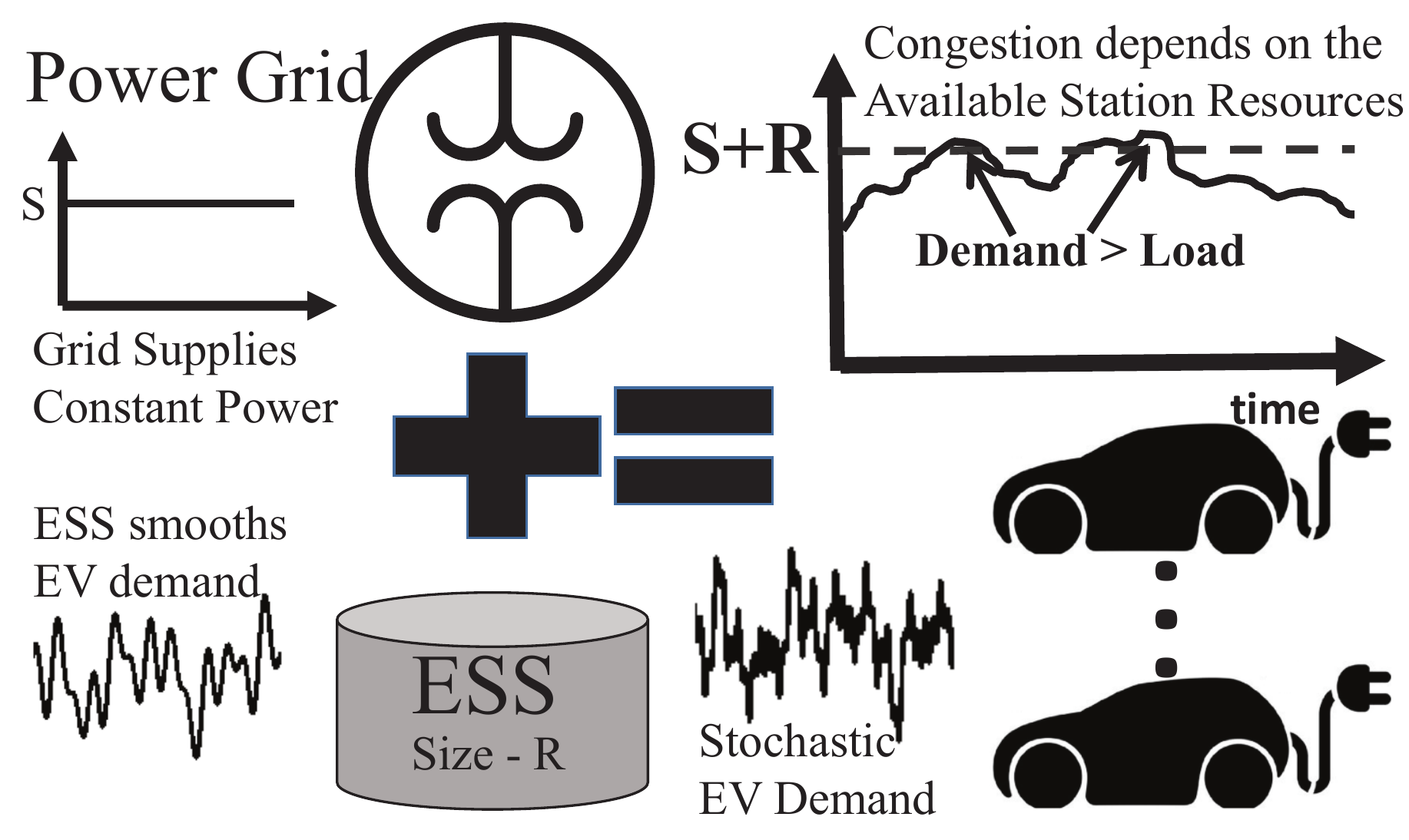}
  \caption{Single charging station architecture.}\label{stationModel}
  \vspace{-15pt}
\end{figure}
\subsection{Previous Work}\label{preWork}
In this study, we build on our previous work \cite{jsac} and \cite{sgc12} that proposed a fast charging station architecture
employing local energy storage and also introduced a stochastic model for customers arrival and service (both from the grid and the local
energy storage device) to assess its performance (depicted in Fig.~\ref{stationModel}). Next, we summarize the dynamics of the single charging station
model introduced in \cite{sgc12}. Each charging station
draws constant power from the power grid, which is expressed by
the number of EVs that can be charged simultaneously, up to
capacity of $S$ vehicles. Similarly a local energy storage is
employed which can accommodate up to $R$ EVs in a
fully charged state. Energy storage will be used to meet spikes in
stochastic power demand that exceed the available grid power
level. Since we always draw constant power from the grid due to
a long term contractual agreement between the grid and the
station, this model also aims to isolate the former from demand
spikes and hence enhances its reliability.

Given this setting, we represent the single charging station model as a two dimensional Markov chain with the following parameters: the arrival of EVs is
modeled as a Poisson process with rate $\lambda$, while charging times both from the grid and/or the storage unit are exponentially distributed
with rate $\mu$. The charging time of the storage unit is also exponentially distributed with rate $\nu$. As mentioned above, EV
charging occurs first from the grid and if that reaches capacity, only then is the storage unit engaged. Finally, if a customer
comes to a charging station when all power/energy resources are used, she could not receive service and consequently will be ``blocked".
 The long-term``blocking probability" represents the QoS performance metric used.

Further, we propose a profit model which prompts station operators to provide a high QoS to customers. In this model, station operators receive
revenue for each served customer ($p_{normal}$) and pay a penalty for each blocked one ($p_B$) and $p_B>p_{normal}$.
 The penalty captures the fact that customer blocking leads to dissatisfaction (part of this penalty could correspond to a voucher given to the customer to alleviate frustration) and degrades the reputation of the station. Combining the revenue and cost components enables station operators to size its capacity to maximize its profit.

This work also builds on \cite{jsac}, that evaluated the above described charging model in a network context using real world traces obtained from the
Seattle public bus system. The obtained results indicated that the spatial distribution of EVs follows a Beta distribution, which is also used in the
present study. Consequently, at given times (e.g., rush hour during weekdays), some regions (e.g., downtown) are busier
than others. Accordingly, stations near a high density area are busier than other stations and unless an EV routing mechanism is in place,
they would fail to meet preset QoS requirements.
\tikzstyle{block} = [rectangle, draw,top color =white , bottom color = processblue!20,
    text width=5em, text centered, rounded corners, minimum height=3em]
\tikzstyle{line} = [draw, -triangle 45]
\begin{figure}[t]
\begin{tikzpicture}[node distance = 2.5cm, auto]
\node [block] (init) {Congestion Level};
\node [block, left of=init] (Grid) {Grid Resources};
\node [block, right of=init] (Policy) {Station Policy};
\node [block, below of=init] (Demand) {EV Demand};
\node [block, right of=Demand] (EVAC) {EVAC (pricing)};
\path [line] (Grid) -- (init);
\path [line] (Policy) -- (init);
\path [line,dashed] (Demand) -- (init);
\path [line] (Policy) -- (EVAC);
\path [line,dashed] (init) -- (EVAC);
\path [line,dashed] (EVAC) -- (Demand);
\end{tikzpicture}
\caption{Interdependence of system components. Allocated grid resources are constant. Station policy determines QoS level. EV demand depends on the time-varying demand and the offered prices.}\label{depend}
\vspace{-15pt}
\end{figure}
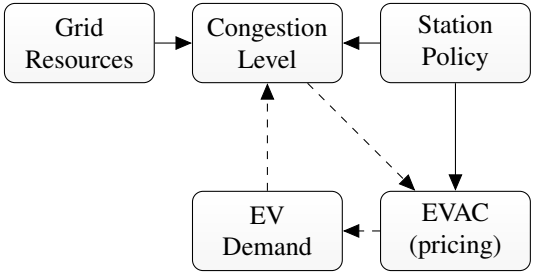

 \begin{figure*}[t]
  \centering
\includegraphics[width=1.2\columnwidth]{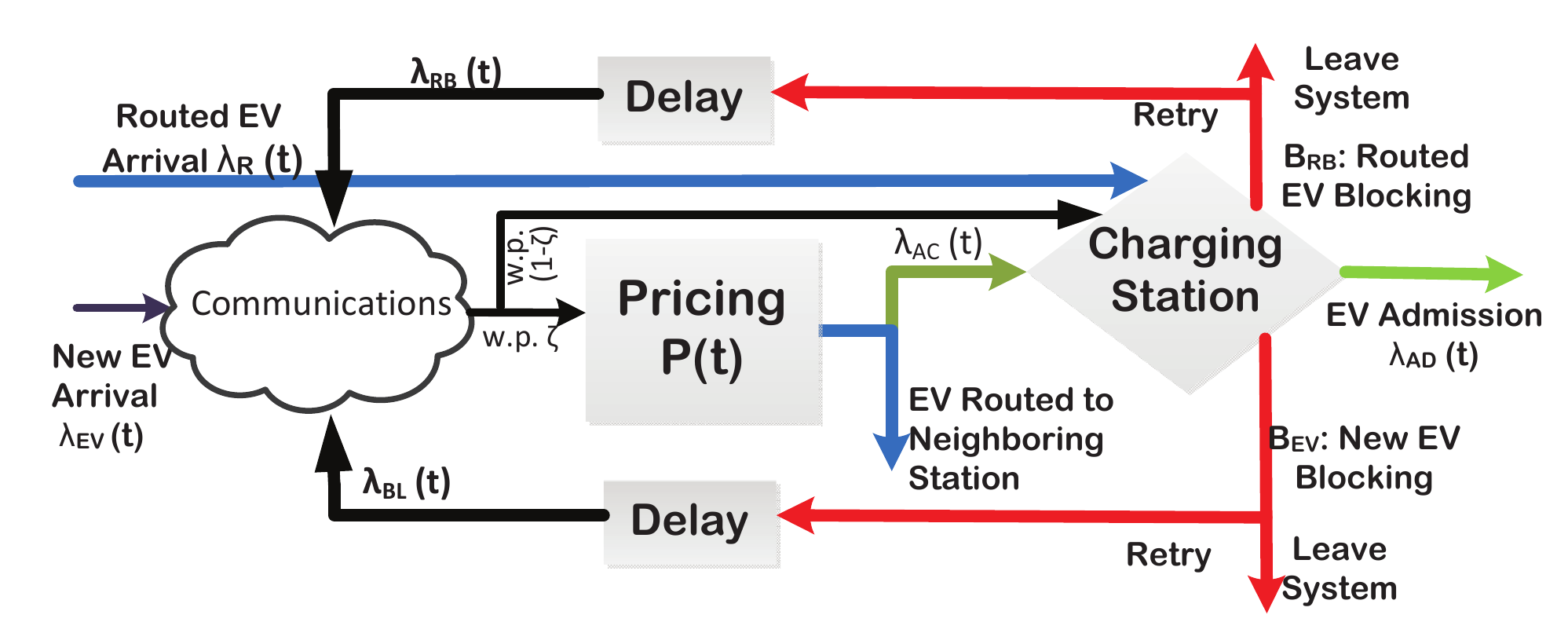}
\caption{Electric Vehicle Admission Control (EVAC) }\label{evac1}
\vspace{-15 pt}
\end{figure*}
\section{Network of Fast Charging Stations}\label{netOfStations}

We present next, the key building blocks of the proposed framework,
starting with Electric Vehicle Admission
Control (EVAC) at a single charging station. Subsequently,
details of the decentralized control scheme for the network of
charging stations are given. An overview of the system components is depicted in Fig.~\ref{depend}.

\subsection{Electric Vehicle Admission Control}

It is assumed that EV drivers act independently and ``selfishly" without
any consideration for the status of the power grid. Thus, uncontrolled
customer demand can easily lead to system failures. Such
situations can be avoided by an admission control mechanism that
offers monetary incentives to customers to alter their choice of a charging station.
The goals of EVAC are: (1) attain QoS targets globally,
(2) balance the EV load among neighboring stations, and (3)
increase the total revenue for the network operator by more
efficient use of grid resources. EVAC relies on a pricing mechanism which adjusts the
price of electricity according to the load on the system. For instance,
during rush hours, station operators can offer relatively lower
prices in a neighboring station to shave the excessive demand at
a busier one.

\subsubsection{System Parameters}

Next, we explain the dynamics of the admission control mechanism employed at a single station level. If there is no congestion it is assumed that each customer chooses to go to the nearest station to receive service. In the EVAC model (depicted in Fig.~\ref{evac1}), the arrival rate for each EV is represented by $\lambda_{\!_{EV}}(t)$. Each customer enters the pricing block, and $\lambda_{\!_{AC}}(t)$ is the proportion of the arrival rate that accepts the offered price. Similarly, $\lambda_{\!_{R}}(t)$ is the arrival rate for routed vehicles from neighboring stations and $\lambda_{\!_{AD}}(t)$ is the arrival rate of admitted EVs. On the other hand some customers may choose to go to busy stations and may not get service. In order to distinguish customer weights, two types of EV blocking probability are defined: (i) $B_{\!_{EV}}$: blocking of EVs who come to the nearest station (ii) $B_{\!_{RB}}$: blocking of a routed EV. Further, when a customer is blocked it is assumed that he will wait some random time and going to retry the pricing block. Overall one third of the customers are assumed to retry and the rest is assumed to leave the system. $\lambda_{\!_{BL}}(t)$ and $\lambda_{\!_{RB}}(t)$ represents the arrival rate for customers who retry the system.

\subsubsection{QoS and Pricing Block}

The blocking probability, which represents the percentage of customers denied service in the long run,
corresponds to the \emph{Quality-of-Service} (QoS) metric. The performance metric is denoted by ($P_{\!_{BT}}$) as the weighted
sum of the two blocking types: ${P_{\!_{BT}}} = {\gamma_1}{B_{\!_{EV}}} + {\gamma_2}{B_{\!_{RB}}} $ where ${\gamma_1} + {\gamma_2} = 1$.
Since blocking an EV that is routed from a neighboring station leads more dissatisfaction, it is assumed that ${\gamma _2} > {\gamma_1}$.
Note that the goal is to provide a charging service within a certain blocking probability.

The pricing block is the key component of EVAC. In general, pricing policies in loss systems can be
classified into three~\cite{Hampshire2009} categories. The first category
uses a \emph{flat pricing} scheme: customers pay by rate or
charged by time of use. The advantage of such policies
is that they are easy to implement and accounting processes are
very simple. On the other hand, such pricing schemes do not
consider the state of the system and hence fail to alleviate
congestions. The second group of policies uses \emph{dynamic
pricing} which monitors the load of the system continuously to
prevent congestions by adjusting the prices. However in the case
of charging stations, implementation of such pricing policies is
rather impractical as it requires the deployment of advanced real time monitoring
and measurement infrastructure.

In our framework, a \emph{myopic policy} that falls in between the first two categories is employed.
Pricing $P(t)$ works as follows. Let
\begin{equation}
{\tilde{\lambda}_{\!_{EV}}}=\lambda_{\!_{EV}}(t)+\lambda_{\!_{BL}}(t)+\lambda_{\!_{RB}}(t),
\end{equation}
then unless the arrival rate $\tilde{\lambda}_{\!_{EV}}(t)$
exceeds a certain threshold $\lambda_{\!_{EV}}^*$ which is
related to the QoS target, the station offers normal price (the
same at each station) $p_{normal}$ that is acceptable by each
customer. During a congestion period (when QoS target is
violated), $\tilde{\lambda}_{\!_{EV}}(t)>\lambda_{\!_{EV}}^*$,
the station owner adjust prices accordingly and offer
\emph{congestion price} $p_{c}>p_{normal}$, so that EVs will
potentially prefer to go to neighboring (also less busy)
stations. Note that the ``peak hour" price is based on
actual customer load on the system.

Next, we define the threshold at which station operator starts to offer congestion prices.
Given the serving capacity of a particular charging station (grid power-$S$, energy storage size-$R$,
and its technology parameter $\nu$), the optimal arrival rate is defined as the maximum arrival threshold
 that the station performance stays within the QoS target $\delta_{max}$. This can be calculated by

\begin{align}\label{optLambda}
{\lambda_{\!_{EV}}^* } = \left\{ {\begin{array}{{l}{l}}
&{\max {\tilde{\lambda}_{\!_{EV}}} }\\
&{\mbox{s.t}.\;{P_{BT}}({\tilde{\lambda}_{\!_{EV}}) \le \delta_{max} }}.
\end{array}} \right.
\end{align}

Note that the QoS target $\delta_{max}$ is specified in the Service Level Agreement (SLA). Next,
 let us define the pricing policy employed in the Pricing Block $P(t)$. As depicted in Fig.~\ref{evac1},
the pricing block $P(t)$ determines the percentage of customers who accept the offered price at time $t$.
That is
\begin{equation}
(\lambda_{\!_{EV}}(t)+\lambda_{\!_{BL}}(t)+\lambda_{\!_{RB}}(t))P(t)=\lambda_{\!_{AC}},
\end{equation}
which yields
\begin{equation}
P(t) = \frac{{{\lambda_{\!_{AC}}}}}{{{\lambda_{\!_{EV}}+\lambda_{\!_{BL}}+\lambda_{\!_{RB}}}}} \le \frac{{\lambda_{\!_{EV}}^*}}{{\tilde{\lambda}_{\!_{EV}}}}.
\end{equation}

This indicates that the percentage of customers who would accept
the offered price is inversely proportional to the load on the
system. Thus, as the arrival rate $\tilde{\lambda}_{\!_{EV}}$
increases, the station operator will increase the prices and less
customers will accept the price.  Note that $D(p_{normal})$=$1$
meaning that all customers accept this price.
It is well-known that demand functions describe the sensitivity of customers to the price changes, that is
\begin{equation}
p(t) = {D^{ - 1}}\left(\dfrac{{\lambda_{\!_{EV}}^*}}{{{\lambda_{\!_{EV}}+\lambda_{\!_{BL}}+\lambda_{\!_{RB}}}}}\right),
\end{equation}
where
$\frac{{\lambda_{\!_{EV}}^*}}{{{\lambda_{\!_{EV}}+\lambda_{\!_{BL}}+\lambda_{\!_{RB}}}}}$
is the load on the system. Various demand functions --for
services providing QoS in congested networks-- have been
proposed in the literature. In this work, we use the demand
function proposed in~\cite{Fishburn}. Then, $p(t)$ the price at
time $t$ becomes:
  \begin{align}\label{congPrice1}
 p(t)  = \left\{ \begin{array}{l}
         p_{normal}, \hfill
          \mbox{if ${\tilde{\lambda}}_{\!_{EV}}(t) \le \lambda_{\!_{EV}}^*$} \\
          \\
        {p_c} = {p_{normal}}\left\{ {1 +\theta \sqrt { - \log \left(\dfrac{{\lambda_{\!_{EV}}^*}}{{{\lambda_{\!_{EV}}+\lambda_{\!_{BL}}+\lambda_{\!_{RB}}}}} \right) }}      \right\}, \\ \hfill \mbox{if ${\tilde{\lambda}}_{\!_{{EV}}}(t) > \lambda_{\!_{EV}}^*$}\end{array} \right. &.
 \end{align}
When congestion occurs in a charging station (${\tilde{\lambda}}_{\!_{{EV}}}(t) > \lambda_{\!_{EV}}^*$) resource pricing
according to the customer load will improve the deterioration of the QoS. Parameter $\theta$ is set by the utility operator and
it will be further explained in the next section.

\subsection{Resource Allocation in a Network of Fast Charging Stations}

As the limitations of the power grid prevent charging stations
from providing more serving capacity, grid operators have to
consider the spatio-temporal aspects of customer demand to optimally
allocate the available grid power. Let us assume that $N$ charging
stations are employed in a serving area and let $l =
{1,2,...,N}$ be the index set and $\mathcal{N}$ be the set of
charging stations. Further, let $S^{max}$ be the maximum portion
of generation capacity that the utility can provide to charging
network. Then, considering the spatial distribution of vehicles
and the discretization assumption at each charging station,
network operator aims to maximize the serving capacity by
optimally allocating the power resources. To this end, network
owner solves a two phase resource allocation problem:

\paragraph{Phase-I}
\begin{align}\label{AllocateS}
& \min_{{S}} &&\sum\limits_{i\in l} P_{BT}^{(i)}(\lambda_i, S_i, R_i)\\
& \text{s.t. }&&\sum\limits_{i\in l} S_i \leq S^{max}\\
& & &\delta_{min} \leq P_{BT}^{(i)}(S_i, R_i, \lambda_i) \leq \delta_{max} \\
& & &S_i \; \in\mathbb{Z^+}, \;\delta_{min},\;\delta_{max}>0\\
& & &R_i\; \mbox{and} \;\lambda_i\; \mbox{are\; given},\; \forall i\in l
\end{align}

Further, given any set of station locations and the EV spatial
distribution, arrival rates for each station
($\lambda_i$) can easily be calculated from a discrete event
simulation~\cite{jsac} (details provided in the
next sections). Also note that the blocking probability has a lower
bound $\delta_{min}$ (typically around 0.01\%) which prevents
capacity planners from over-provisioning serving capacity.
However, in reality, the result of the optimization problem
($\boldsymbol{S}$=($S_{1}$,...,$S_{N}$)) is constrained by the power grid operation
limitations. Hence, another constraint $S_i \leq S^{limit}$
that dictates the upper limit for the grid power that can be
drawn, is introduced. Then, if $\sum\limits_{i\in l} S_i <
S^{max}$, meaning that one or more stations cannot get optimal
$S$ due to power grid limitations, the excessive power,
\begin{equation}
S^{excess}=S^{max}-\sum\limits_{i\in l} S_i ,
\end{equation}
is allocated among the neighboring stations inversely
proportional to the physical distance between stations.

Let $\Psi\subset \mathcal{N}$ be the subset charging stations
which have excessive power. Suppose that station $j\in\Psi$ has
excessive power $S_j^{excess}$. Then in Phase-II a neighboring
station-$k\in \mathcal{N}\setminus\Psi$ gets an extra
$S_{kj}^{extra}$ amount of power from station-$j$, that is:

\paragraph{Phase-II}
\begin{equation}\label{excessive}
S_{kj}^{extra} = \ceil[\big]{S_j^{excess}{\omega_k}},\mbox{where}\;{\omega_k} = \dfrac{1}{{d_{kj}^2}}\left(\sum\limits_{z \in N \setminus \Psi } {{{ {\frac{1}{{d_{zk}^2}}} }}}\right)^{ - 1} ,
 \end{equation}
where $d_{zj}$ is the physical distance between stations z and j. A simple example is discussed next, to better clarify this point. Assume that four stations serve a neighborhood, and station 1 has excessive power of six units. The remaining three stations do not have any excessive power and they will share six units of power according to \eqref{excessive}. Let us assume that $d_{21}$=$1$, $d_{31}$=$2$, and $d_{41}$=$2$ units of distance. Then, the excess power would be allocated as follows: $S_{21}^{extra}$=$4$, $S_{31}^{extra}$=$1$, and $S_{41}^{extra}$=$1$. Notice that the power allocation is inversely proportional to the inverse of the square of the physical distance. The main reason is that in the event of congestion, offering enough incentives to customers in order to convince them to drive much further away from their current location is a rather challenging task. Thus, since the closest station (in this case station 2) will serve most of station 1's customers, it gets a bigger portion.
\begin{table}[t]
  \centering
  \caption{Systems Parameters}
    \begin{tabular}{p{1.0cm}  p{\columnwidth-2cm} }
    \toprule
    \multicolumn{2}{c}{} \\

    Parameter & Definition \\
    \midrule
    $\lambda_{\!_{EV}}(t)$ & EV arrival rate (charge request) to the nearest station.\\
    $\lambda_{\!_{AC}}(t)$ & EV arrival rate who accepted the price at the nearest station.\\
    $\lambda_{\!_{R}}(t)$ & Arrival rate of routed (from neighboring stations) customers.\\
    $\lambda_{\!_{BL}}(t)$ & Arrival rate of blocked customers--specific percentage of EVs retry the system.\\
    $\lambda_{\!_{RB}}(t)$ & Arrival rate of blocked routed customers --specific percentage of EVs retry the system.\\
    ${B_{\!_{EV}}} $ & Blocking probability of customers (who pick the nearest station).\\
    ${B_{\!_{RB}}} $ & Blocking probability of the routed customers.\\
    $p_{normal}$ & Normal price accepted by every driver. The same price at each station.\\
    $p_{congestion}$ & Congestion price depends on the station load. \\
    $p_{\!_B}$ & Penalty paid for blocking a customer. Alleviates customer frustration.\\
    $\delta$ & QoS target may be different at each station. \\
    $\lambda^*$ &Maximum allowable arrival rate to support QoS calculated by~\eqref{optLambda}.\\
    $c_{k_{inctv}}$ & Minimum amount of desired savings to accept routing. \\
    $c_{k_{near}}$ & Cost to drive to nearest station. \\
    $d_k$ & Distance vector for each customer k.\\
    $\theta$ & Vector of price tuning parameter for each station set by station operator. \\
    $S$, $R$, $\mu$ & Single station parameters. Denotes grid power, energy storage size, and EV charging rate respectively.\\
    $S^{limit}$ & Maximum power that can be drawn by a single station.\\
    \bottomrule
    \end{tabular}%
  \label{systemParam}%
  \vspace{-10pt}
\end{table}%

\subsection{Decentralized Control}

The motivation for proposing a decentralized control mechanism is to analyze how to serve a large population of EVs under QoS guarantees without requiring significant upgrades to grid resources. For this purpose, we employ a game theoretic framework and propose a decentralized control strategy to balance the EV load among neighboring charging stations. The players of the game are the network operator and the EV drivers. The objective of the ``self-interested" EV drivers is to receive service within QoS guarantees and to minimize the cost of service, which is proportional to the charging price and the distance to the charging facility. On the other hand, the network operator aims to maximize the number of served customers by incentivizing EV drivers to cooperate and route customer demand to neighboring stations with relatively lower utilization, if necessary.

As depicted in Fig.~\ref{evac1}, the proposed framework relies extensively on the capability of both parties to communicate. It is assumed that each charging facility can communicate with a central unit that can forward pricing signals to mobile EVs, and influence their preferences. On the other hand the EVs respond to these signals and play a best response strategy. In this paper, we assume that there exists a high fidelity communications infrastructure. However, as part of future work, we plan to explore the candidate communication technologies and quantify the effects of their quality of service to the charging network performance.

\subsubsection{Game Formulation}
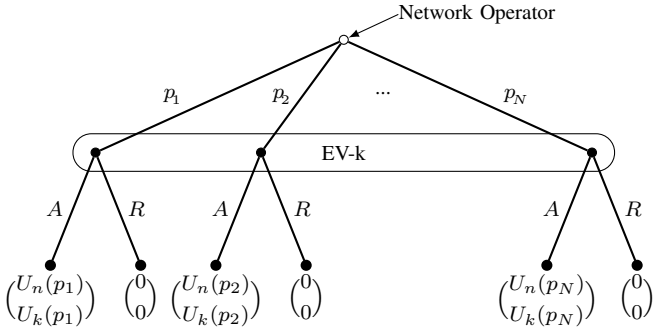
\begin{figure}[t]
\centering
\begin{tikzpicture}[font=\footnotesize,edge from parent/.style={draw,thick}]
\tikzstyle{solid node}=[circle,draw,inner sep=1.2,fill=black];
\tikzstyle{hollow node}=[circle,draw,inner sep=1.2];
\tikzstyle{level 1}=[level distance=15mm,sibling distance=22mm]
\tikzstyle{level 2}=[level distance=15mm,sibling distance=12mm]
\node(0)[hollow node]{}
child{node[solid node]{}
child{node[solid node]{}edge from parent node[left]{$A$}}
child{node[solid node]{}edge from parent node[right]{$R$}}
edge from parent node[left,xshift=-10]{$p_{\!_1}$}
}
child{node[solid node]{}
child{node[solid node]{}edge from parent node[left]{$A$}}
child{node[solid node]{}edge from parent node[right]{$R$}}
edge from parent node[left,xshift=0]{$p_{\!_2}$}
}
child{node[draw=none]{} edge from parent[draw=none] node[align=center]{...}}
child{node[solid node]{}
child{node[solid node]{}edge from parent node[left]{$A$}}
child{node[solid node]{}edge from parent node[right]{$R$}}
edge from parent node[right,xshift=10]{$p_{\!_N}$}
};
\draw[rounded corners=7]($(0-1)+(-.3,.25)$)rectangle($(0-4)+(.3,-.25)$);
\draw[<-,>=latex](0)--(25:8mm)node[inner sep=0,right]{Network Operator};
\node at($.5*(0-1)+.5*(0-4)$){EV-k};
\node(payoff)[below]at(0-1-1){$\displaystyle\binom{U_n(p_1)}{U_k(p_1)}$};
\node[below]at(0-1-2){$\displaystyle\binom{0}{0}$};
\node[below]at(0-2-1){$\displaystyle\binom{U_n(p_2)}{U_k(p_2)}$};
\node[below]at(0-2-2){$\displaystyle\binom{0}{0}$};
\node[below]at(0-4-1){$\displaystyle\binom{U_n(p_N)}{U_k(p_N)}$};
\node[below]at(0-4-2){$\displaystyle\binom{0}{0}$};
\end{tikzpicture}

\caption{Stackelberg game in its extensive
form.}\label{extensiveForm}
\vspace{-15pt}
\end{figure}

 Customers have access to charging station location information and the corresponding pricing offers via the central communication unit. As outlined above, the role of the network owner is to communicate with each charging station and forward their price signals; the ultimate goal is to attract customers to drive to a more distant station to balance the arrival rate at each station. Thus, the charging network owner acts as a \emph{leader} who can commit to a strategy before \emph{followers} (customers) can pick their strategies. In this respect, a Stackelberg game (leader-follower game) is employed to model this system. First, we define the game in its strategic form $\Gamma=\{ \{\mathcal{N}\cup \mathcal{K}\}, \{  \boldsymbol{p}_{k} \} , \{X_{k\in\mathcal{K}}\}, \{U_{n\in\mathcal{N}}\}, \{U_{k\in\mathcal{K}}\}\}$, where $\mathcal{N}$=$\{1,\cdots, N\}$ is the set of charging stations and $\mathcal{K}$=$\{1,\cdots, K\}$ the set of EVs that require charging at a given time. The strategies of each set of players are as follows. $\{\boldsymbol{p}_k\}$ denotes a $1\times N$ price vector offered by the charging network operator. Note that to set the pricing vector station operator also need to set $\{\boldsymbol{\theta}\}$ for each station. The rationale behind setting this parameter stems from the fact that in order to convince people to change their electricity consumption during rush hours, or in our case where to charge their EVs, the ratio of peak-to-off-peak hour should be some positive number that is greater than one \cite{brattle}. The parameter $\theta$ is set to maintain such incentives for customers. The EVs strategy, $X_{k}$, is to pick a charging station from $\mathcal{N}$. $U_{n}$ and $U_{k}$ are the utility functions that represent the payoffs for the both players. The game in its extensive form is given in Fig.~\ref{extensiveForm}. Next, we describe the components of the game in detail.

\subsubsection{EV Customers - Followers}
In the absence of any incentives and under the assumption of rationality, each customer will just choose to go to its nearest charging station. Through the use of two-way communications, the network operator aims to route EVs to less busy stations (if necessary) by offering relatively lower prices, but still avoid penalties for poor QoS (blocked EVs). The EVs
can either \emph{Accept} to go to a less busy charging station or \emph{Reject} and go to the nearest one. The decision will be driven by the price responsiveness of the customers. As presented in~\cite{faruqui2011will}, in order to motivate EV owners to change their charging behavior, there should be enough financial incentives. Considering this, the EV strategy is given below.
 \begin{align}
 EV_{Strategy} = \left\{
 \begin{array}{r r r r}
    \mbox{Accept,}& \text{if $c_{\!_{k_{near}}}-c_{\!_{k_{desrd}}}\geq c_{\!_{k_{inctv}}}$ } \mbox{and}\\
    & P_{BT}^{desired}\leq \delta_{max}  \\
    \mbox{Reject,} & \text{otherwise.} \\
  \end{array} \right.
\end{align}

Suppose that customer $k$, $k\in\mathcal{K}$  wants to get a
service from station $n$, $n\in\mathcal{N}$. The cost of
choosing the nearest station equals to $c_{\!_{k_{near}}}$=$p_n$
+ $c\left( d_{nk}^2 \right)$, where $p_{n}$ is the price paid
for the service at station $n$, and $c \left(d_{nk}^2\right)$ is
the cost of driving to station $n$ for EV $k$. Note that $p_{n}$
equals $p_{normal}$  when there is no congestion
(${\tilde{\lambda}_{\!_{{EV}}}(t)} \le \lambda_{\!_{{EV}}}^*$)
at the nearest station $n$, and to the congestion price $p_c$
otherwise. Similarly, when the network operator desires to route
a vehicle to another station, it costs $c_{\!_{k_{desrd}}}$. As
described above there should be strong enough incentives for drivers to
drive extra miles. Hence, to project customer behavior into our
formulation, we assume that EV $k$ will only agree to go to the
desired station if it gains at least $c_{\!_{k_{inctv}}}$ (e.g.,
considering current electricity prices there should be a
reasonable level of savings e.g., $10-15\%$ savings).

The pay-off of a single EV is represented by a utility function that has the following three parameters: the QoS metric $\boldsymbol{P}_{\!_{BT}}$, the price at each station $\boldsymbol{p}_k$, and the physical distance to each station $\boldsymbol{d}_k$. Note that each component of the utility function corresponds to a $1\times N$ row vector, where the position of each element is associated with the corresponding charging station parameter (QoS, service price, and distance). Also the pricing scheme proposed in section~\label{PricingBlock} is employed at each station. Then, the utility function for EV-$k$ becomes:
\begin{equation}
{U_k}(\boldsymbol{{P}}_{\!_{BT}},\boldsymbol{d}_k,\boldsymbol{p}_k) = h(\boldsymbol{{P}}_{\!_{BT}}) \boldsymbol{e}_n \left\{ \boldsymbol{p}_k + {c_k}(\boldsymbol{d}_k) + {f_k}(\boldsymbol{d}_k) \right\},
\end{equation}
where $h({\boldsymbol{P}_{\!_{{BT}}}})=e^{\xi ({\boldsymbol{P}_{BT}} - \delta )}$. The $h( \cdot )$ function denotes the disutility of experiencing high blocking probability \cite{lau1998mobility}, $ \delta$ is the QoS target and $\xi$ is a constant which reflects how much the customer values the service (e.g., higher value for when in urgent need of charging service). This component is used because blocked customers may retry to get service and high blocking leads further disutility. These components capture the dissatisfaction introduced
by a high blocking probability. Note that,
\begin{align}
\xi = \left\{ {\begin{array}{*{20}{c l}}
0,&\lambda \leq \lambda^{*} \\
\xi \in \mathbb{R^{+}}, &\text{otherwise.}
\end{array}} \right.
\end{align}
Customer $k$ picks one station $n\in\mathcal{N}$ and $\boldsymbol{e}_n$
represents a column vector comprising of all zeros except for
the $n^{th}$ position which is $1$. Similarly, $\boldsymbol{p}_k$ is
the price vector offered by the charging network operator. The
next term in the utility function reflects the cost of driving
to a charging station which is a function of distance to each
station from the current location of EV $k$. Finally,
${f_k}(\boldsymbol{d}_k)$ is related with the dissatisfaction of EV $k$
when it selects to go to a more distant station. Even though the
total cost is lowered, some level of dissatisfaction occurs due
to spending extra time to drive the extra miles (e.g., time cost
etc.).

Thus, the optimization problem of the EVs is a mapping $ \mathbb{R}  \to \left\{ {{\rm{reject, accept}}} \right\}$ that gives rise to a vector $\boldsymbol{e}_n $,
where ($Accept$) sets the $n^{th}$ position to $1$ while all other remain $0$ ($Reject$). Note that the objective is to minimize the overall cost. This can be expressed as follows:

\begin{align}\label{customer_util}
&\mathop {\mathop {\arg \min}\limits_{n}}&&  h(\boldsymbol{P}_{\!_{BT}})\boldsymbol{e}_n \left\{ \boldsymbol{p}_k + {c_k}(\boldsymbol{d}_k) + {f_k}(\boldsymbol{d}_k) \right\}&\\
&\text{s.t.}&& n\in\mathcal{N} &
\end{align}

where EV $k$ chooses station $n$ and pays $p_n$ amount of money.
The second component is the cost of driving:
${c_k}(\boldsymbol{e}_n\boldsymbol{d}_k)$ = $p_{\!_{dr}}d_{kn}^2$ (cost of
driving one unit of distance times distance to station $n$).
Similar approaches are proposed in the classic analysis of
Hotelling where customers aim to minimize their travel costs to
get service~\cite{hotelling}. The last parameter reflects the
dissatisfaction of customers to drive extra miles. A linear
dissatisfaction model is employed
$f_k({\boldsymbol{e}_n\boldsymbol{d}})$=$p_{dis}(d_{kn}-d_{nearest})$ where
$p_{dis}$ is the cost of driving one unit of distance, and
$(d_{kn}-d_{nearest})$ is the total amount of extra miles
traveled. This component captures the real world
behavior of customers, since it is unlikely that one would drive
a significant number of extra miles for limited savings. Finally, the
above optimization problem is subject to the EV having enough stored charge
to drive to station $n$, which we assume that this
is the case (about $1$-$2$ kWh of stored energy).
\begin{figure}
                \centering
                \includegraphics[width=\columnwidth]{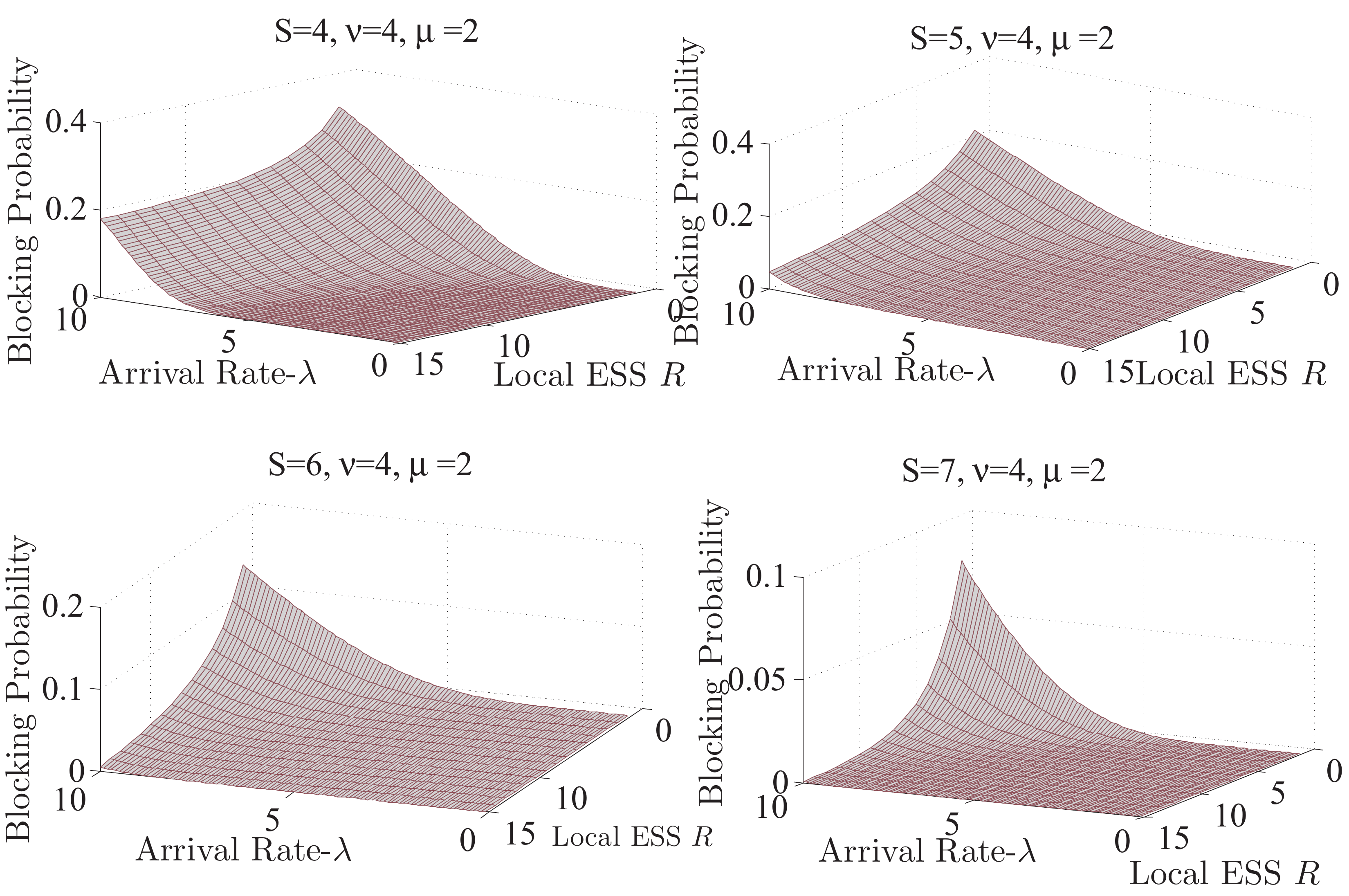}
  \caption{Blocking probability (\% of not served customers) computation for different station parameters.}
  \label{singleBlocking}
  \vspace{-10pt}
  \end{figure}

\begin{figure*}[t]
        \centering
        \begin{subfigure}[b]{0.35\textwidth}
                \centering
                \includegraphics[width=\columnwidth]{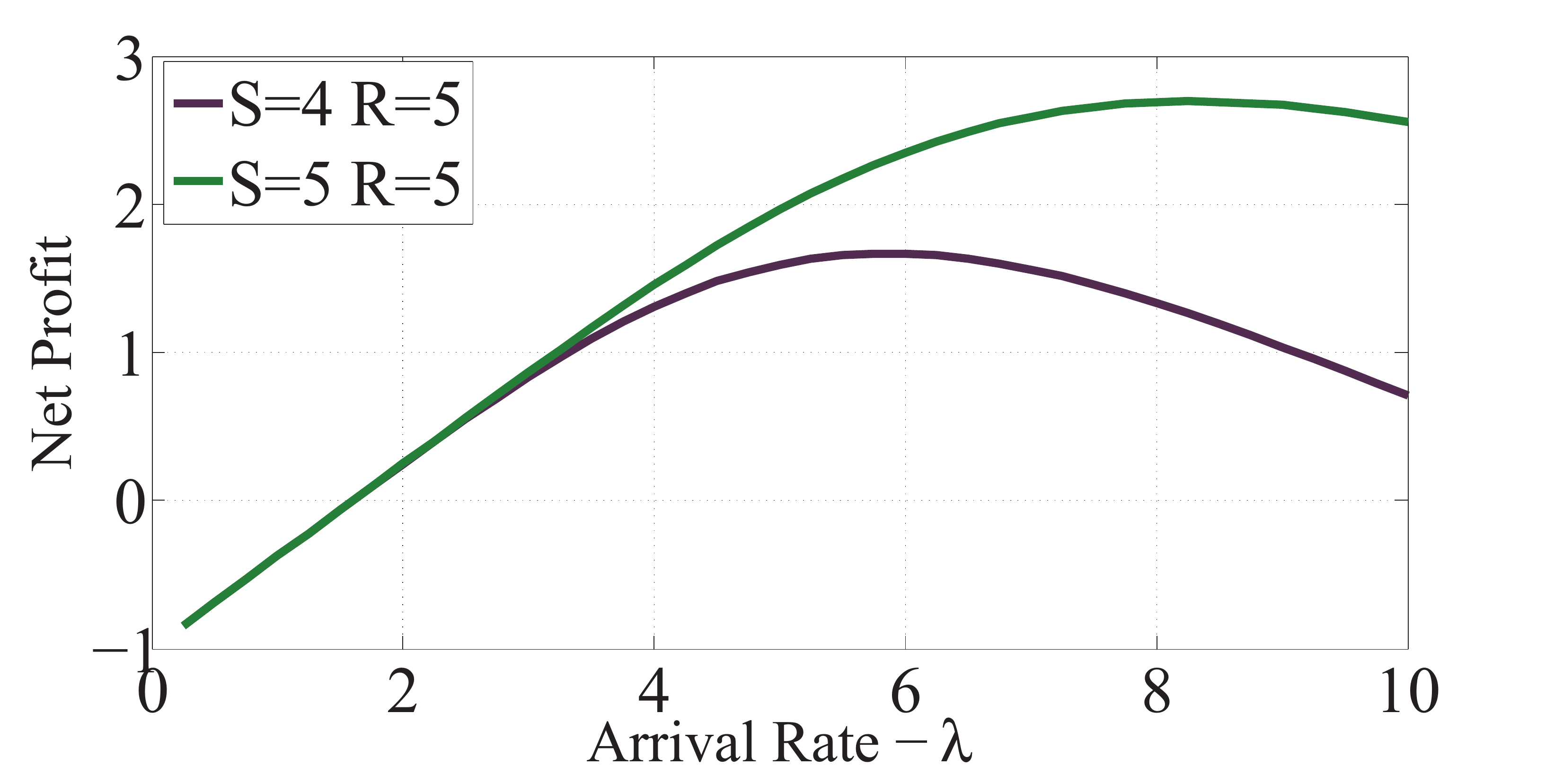}
  \caption{Numerical evaluation-I.}\label{singleProfit}
          \end{subfigure}
        \begin{subfigure}[b]{0.33\textwidth}
                \centering
                \includegraphics[width=\columnwidth]{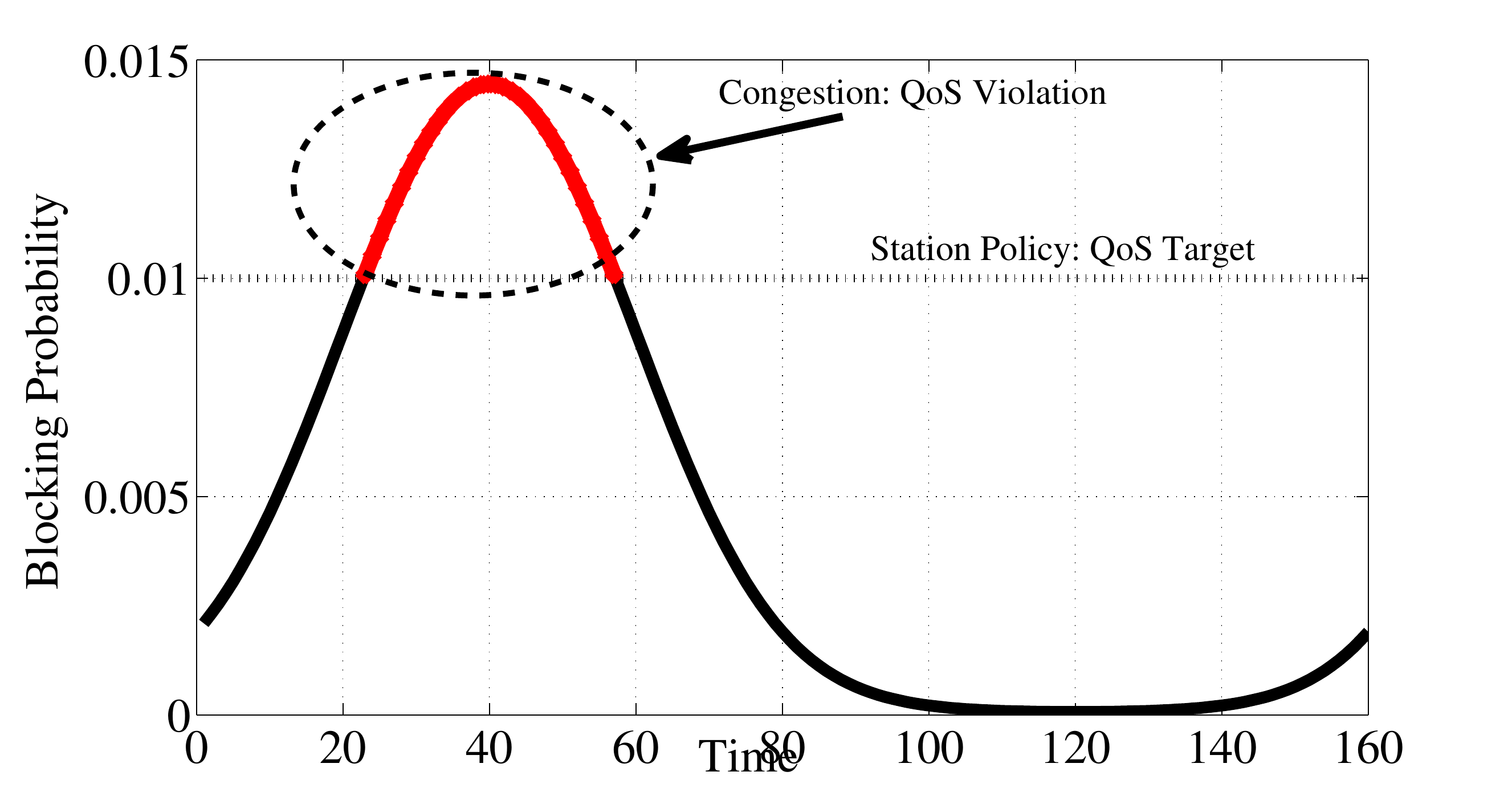}
  \caption{Illustration of single charging station congestion, $\lambda$=$5$+$1.75\sin (2\pi t/80)$. }\label{sineWave}
        \end{subfigure}
 \begin{subfigure}[b]{0.29\textwidth}
                \centering
                \includegraphics[width=\columnwidth]{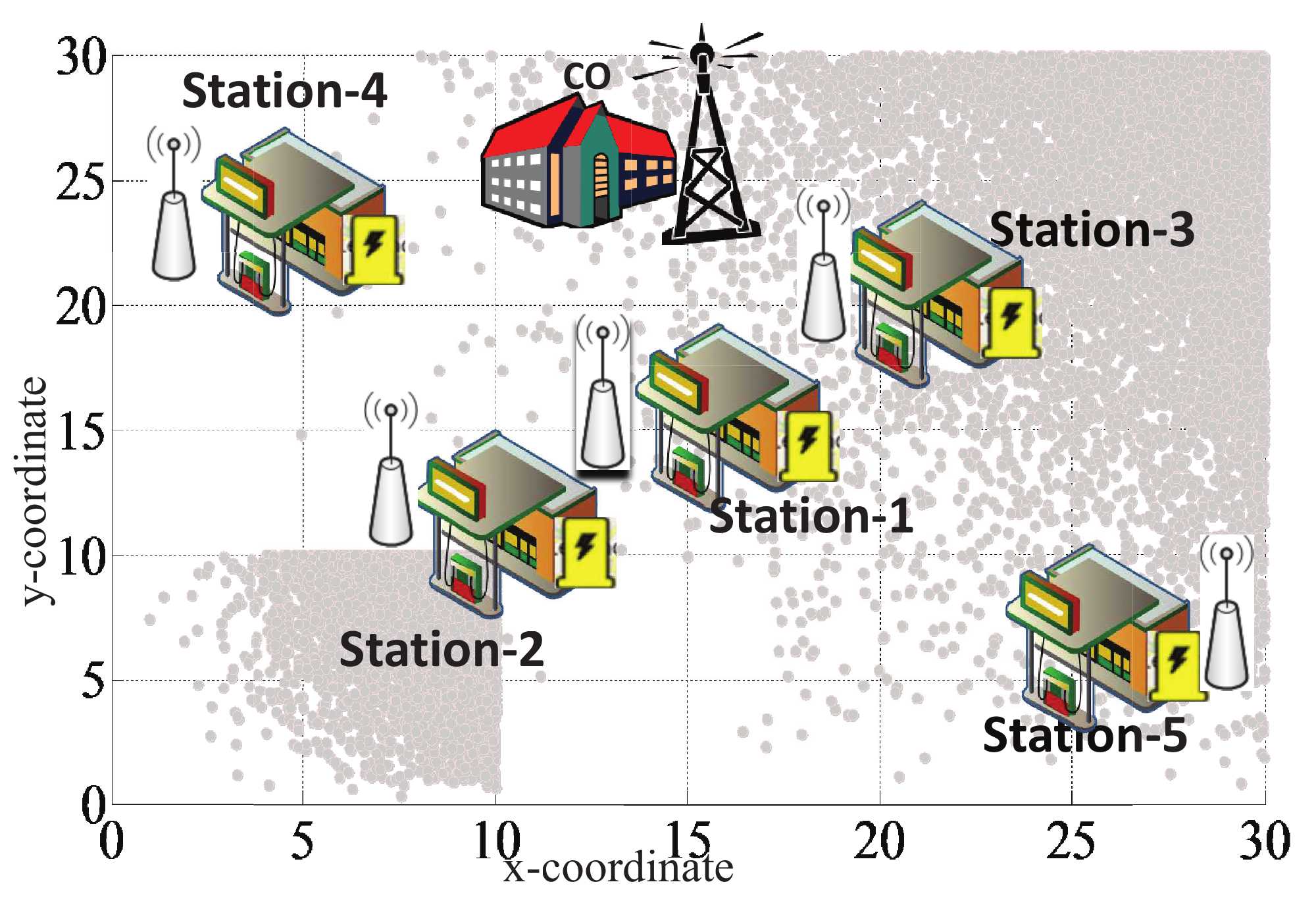}
                \caption{Charging network layout.}
                \label{mapLayout}
  \end{subfigure}
        \caption{Evaluation of single charging station dynamics.}\label{results2}
        \vspace{-15 pt}
\end{figure*}
 \begin{figure*}[t]
        \centering
                \begin{subfigure}[b]{0.32\textwidth}
                \centering
                \includegraphics[width=\columnwidth]{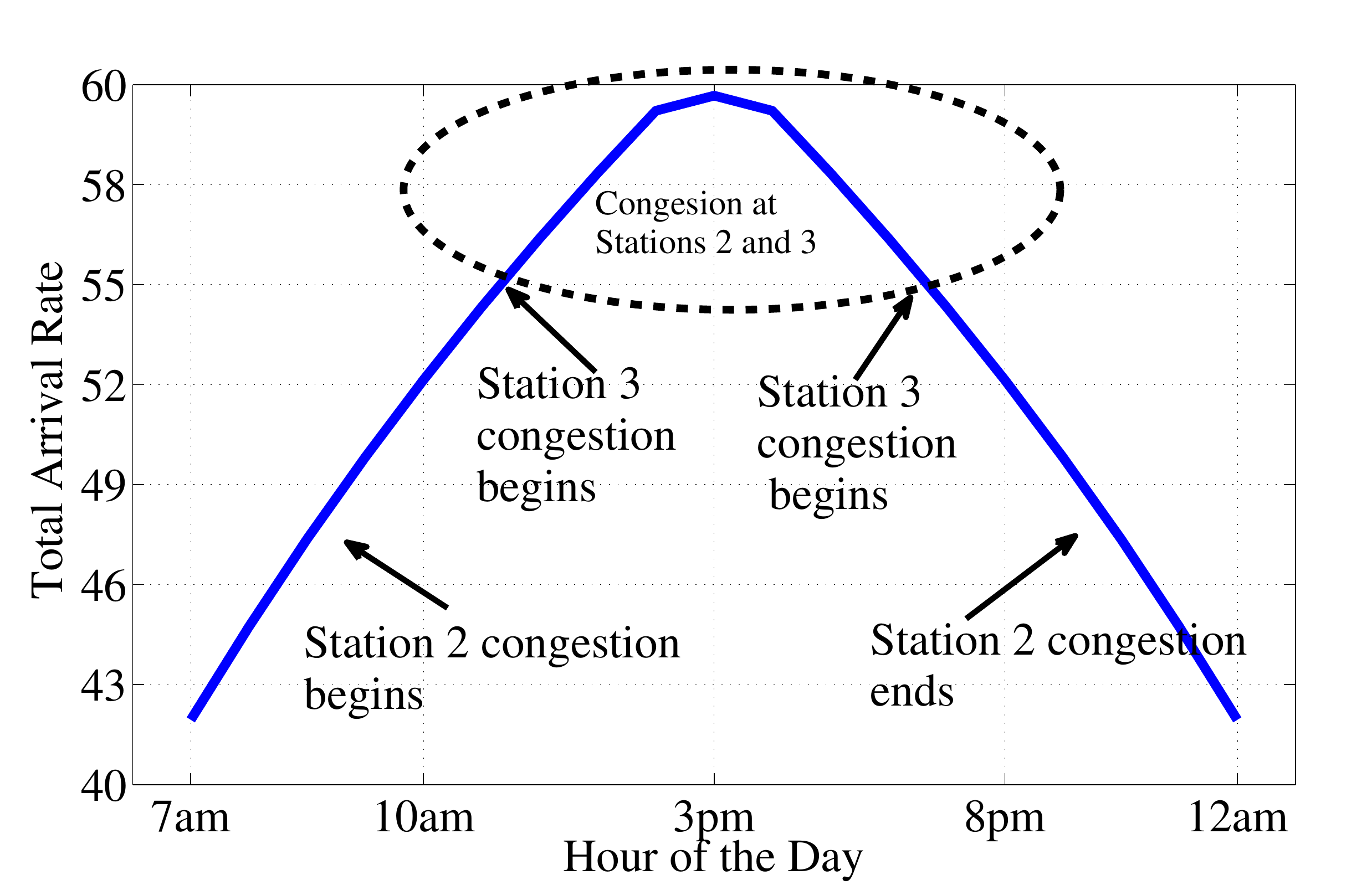}
                \caption{Daily customer demand.}
                \label{custDemand}

        \end{subfigure}
             \begin{subfigure}[b]{0.33\textwidth}
                \centering
                \includegraphics[width=\columnwidth]{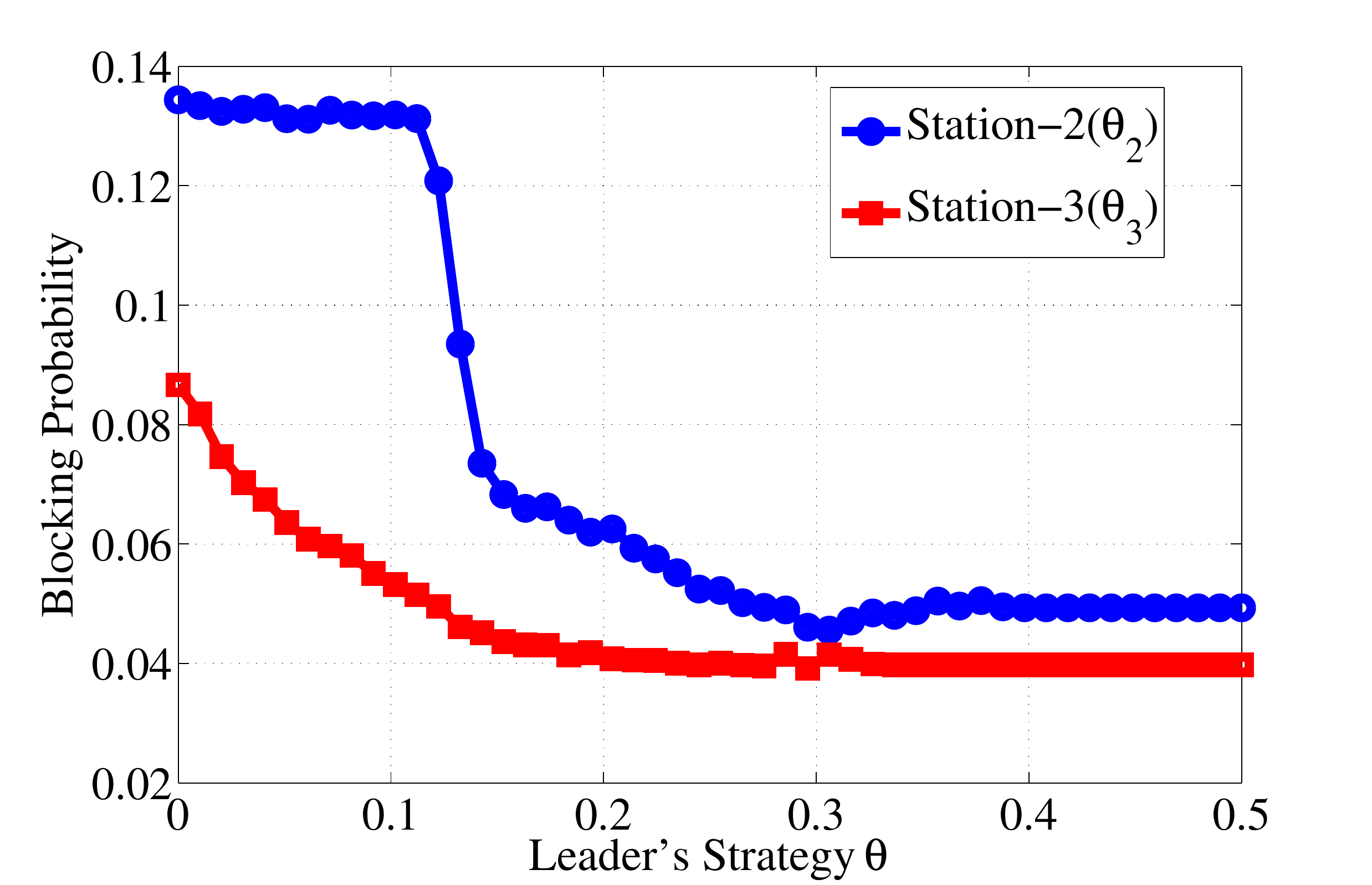}
                \caption{Leader's strategy.}\label{leaderStrategy}
        \end{subfigure}
\begin{subfigure}[b]{0.33\textwidth}
                \centering
                \includegraphics[width=\columnwidth]{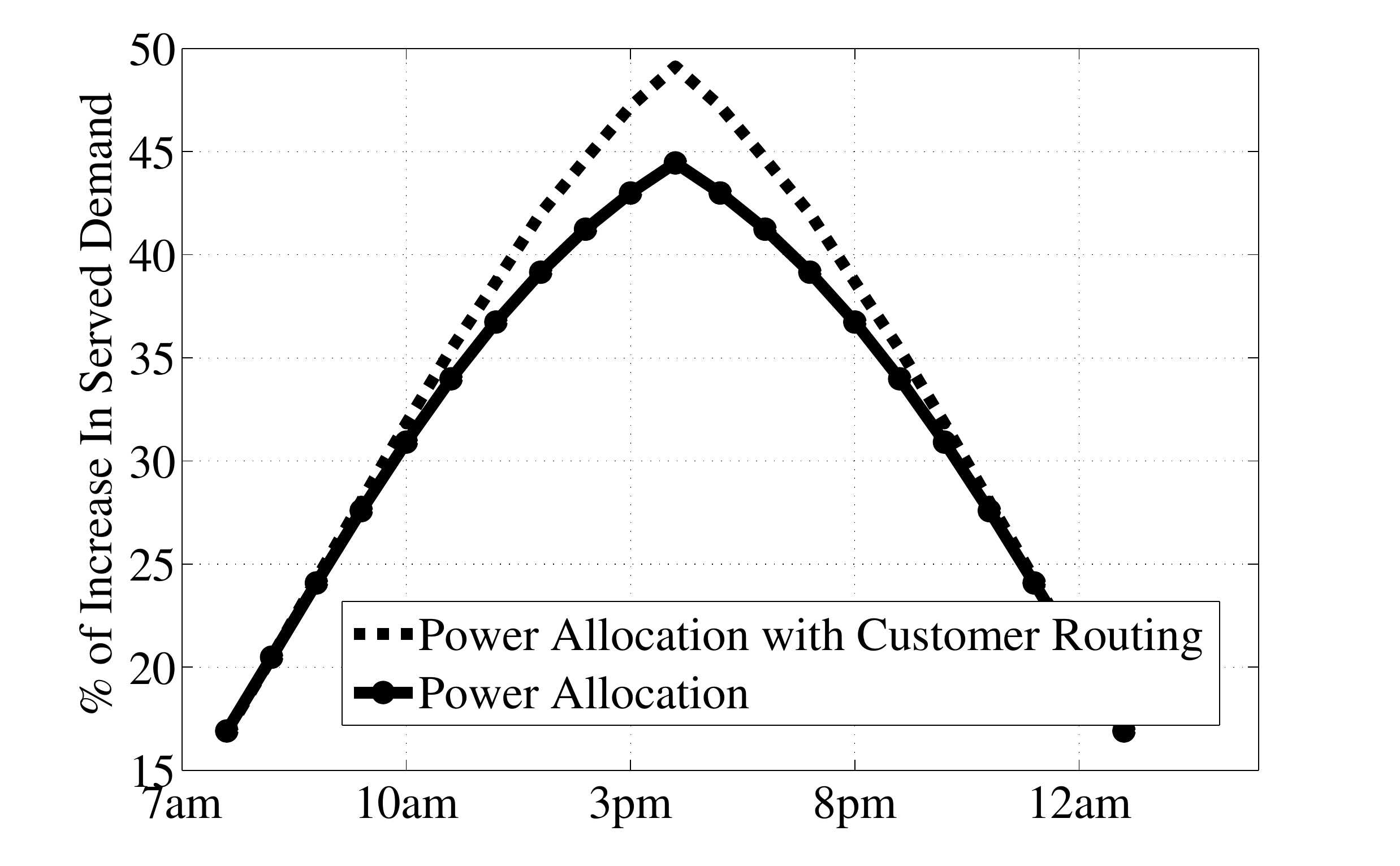}
                \caption{\% of increase in served EVs.}
                \label{IncreasedServings}
        \end{subfigure}%

        \caption{Numerical evaluation-II.}\label{results1}
        \vspace{-15pt}
\end{figure*}
 \begin{figure*}[t]
        \centering
        \begin{subfigure}[b]{0.31\textwidth}
                \centering
                \includegraphics[width=\columnwidth]{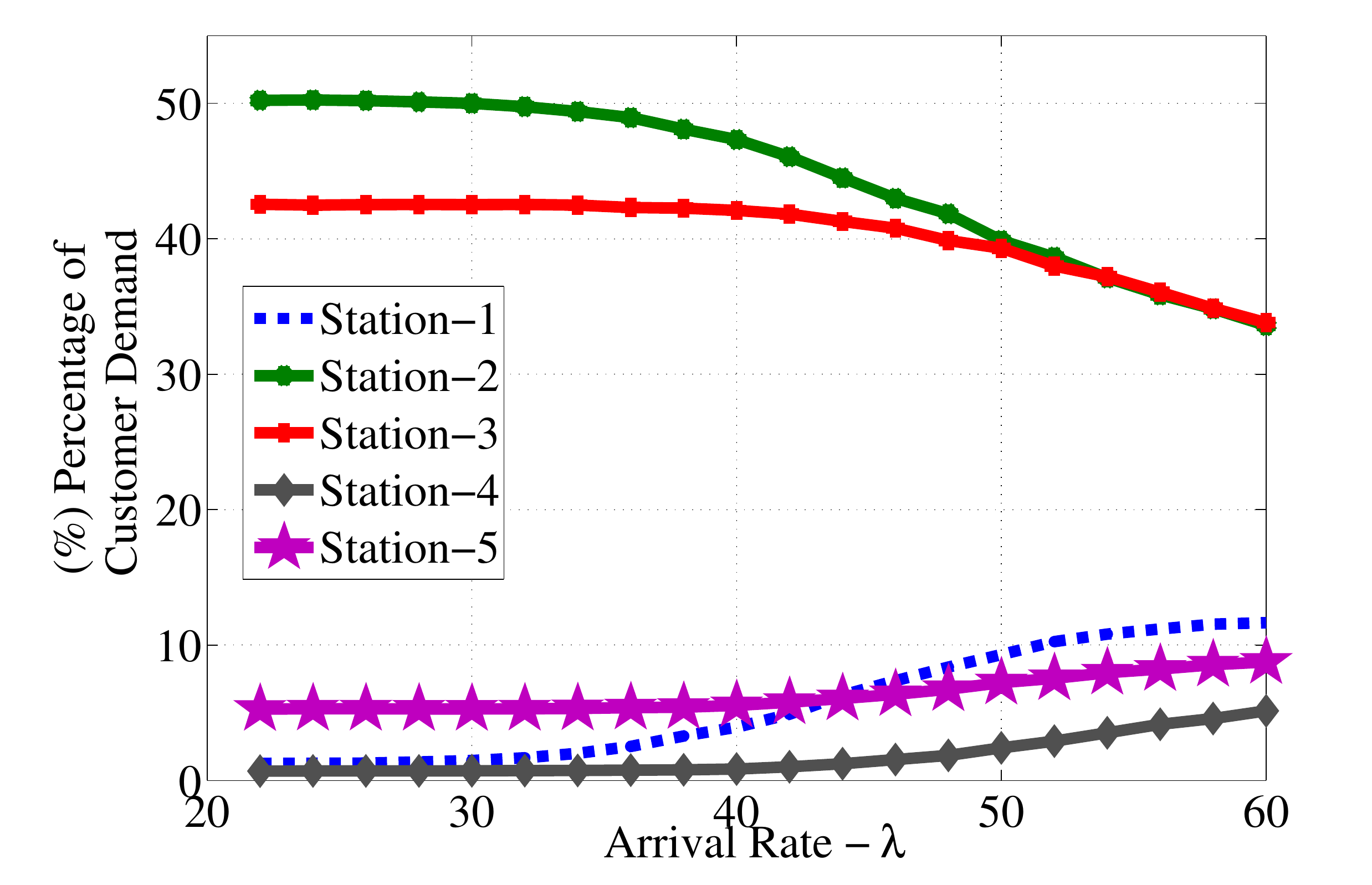}
  \caption{Load balancing in a network of fast charging stations.}
  \label{loadBalancing}
  \end{subfigure}
        \begin{subfigure}[b]{0.32\textwidth}
                \centering
                \includegraphics[width=\columnwidth]{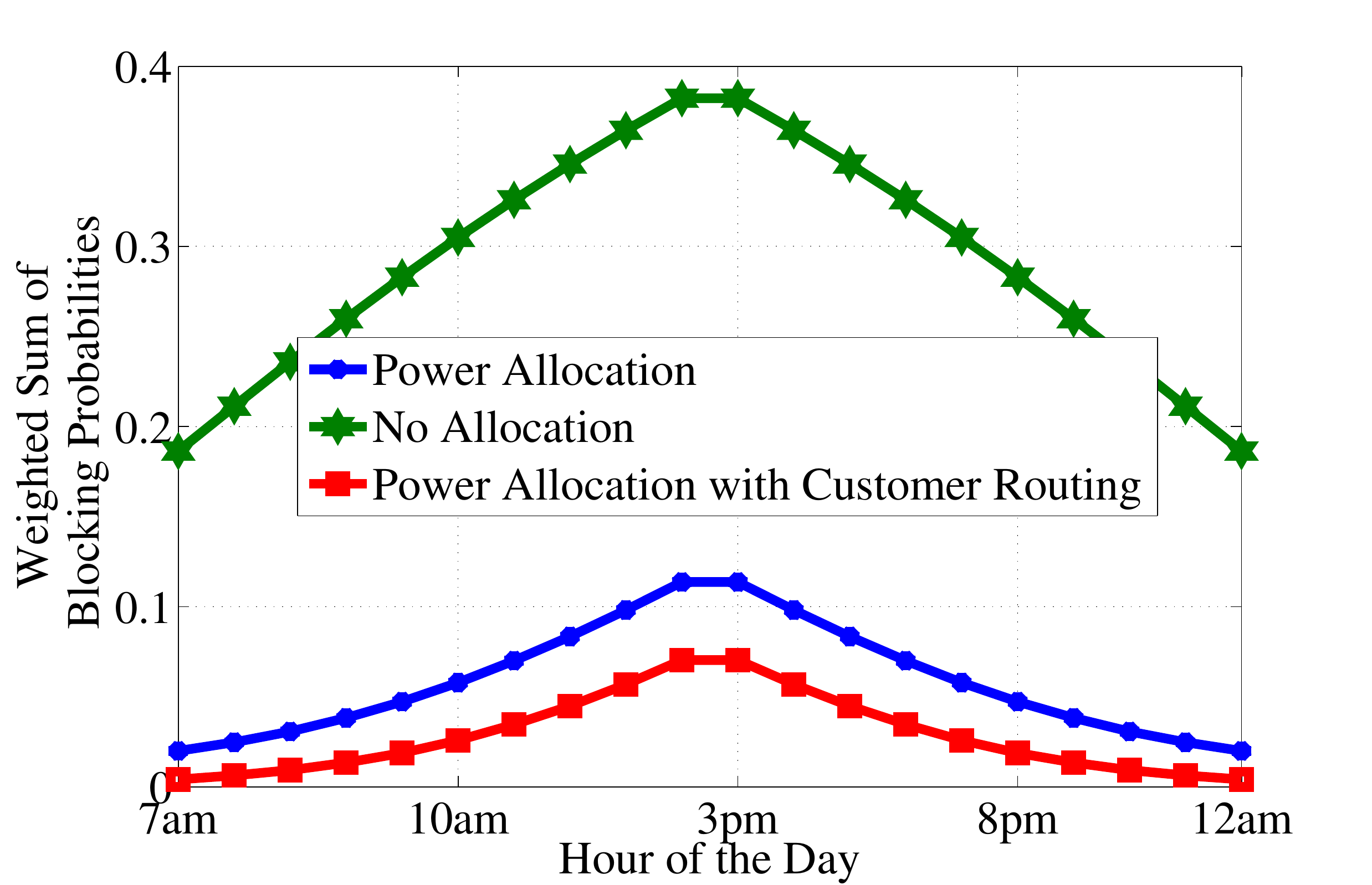}
  \caption{Weighted sum of blocking probability for all stations.}\label{station2}
          \end{subfigure}
        \begin{subfigure}[b]{0.33\textwidth}
                \centering
                \includegraphics[width=\columnwidth]{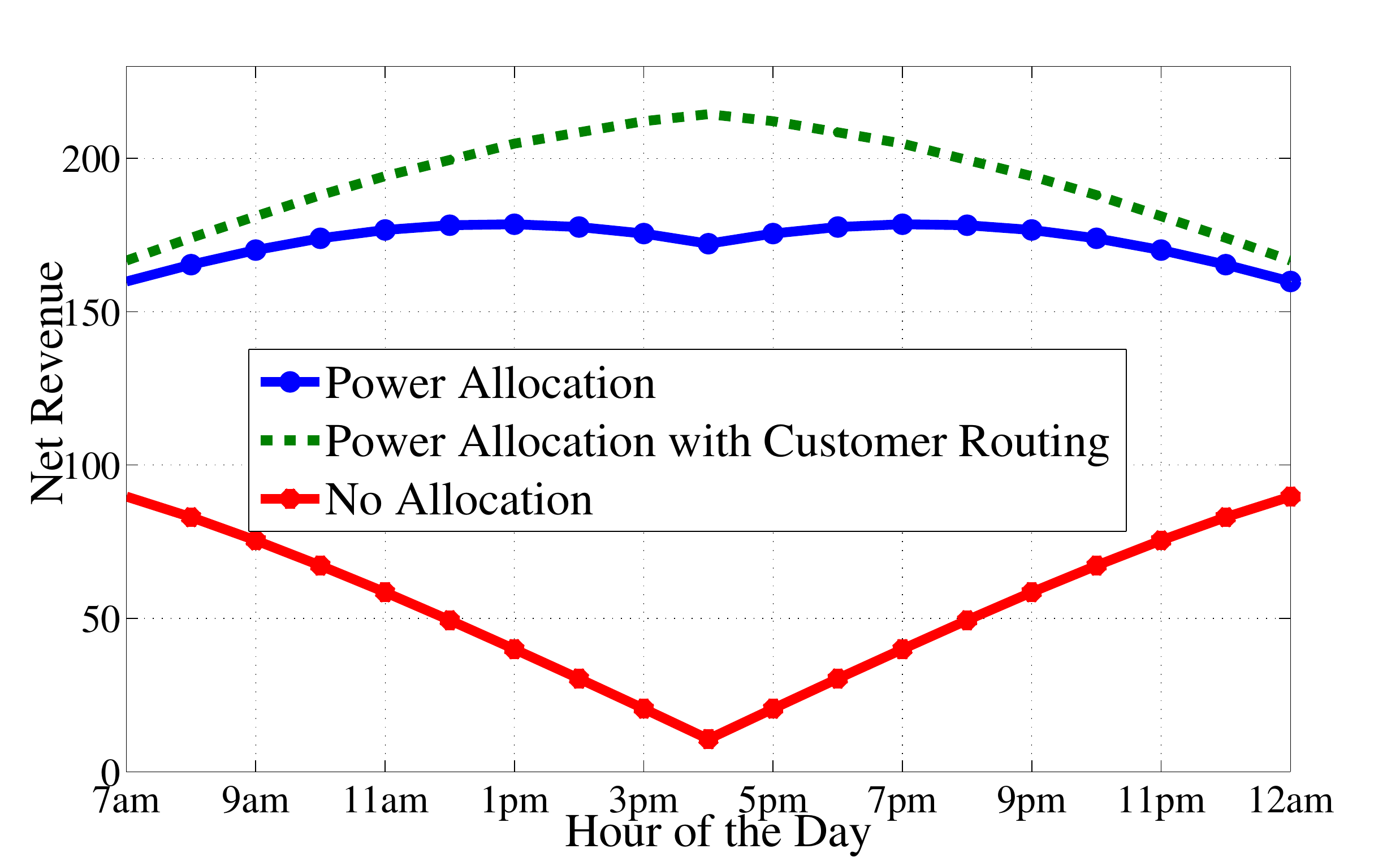}
  \caption{Network-wide (five stations) net revenue comparison. }\label{netRev}
        \end{subfigure}

        \caption{Numerical evaluation-III.}\label{results2}
       \vspace{-20 pt}
\end{figure*}

\subsubsection{Network Operator - Leader}

As described in the previous section, network operator has the first move advantage to set the price and he tries to influence the customer
preferences. In the Stackelberg game, the network operator acts as the \emph{leader} who aims to coordinate non-cooperative EV charging behavior
so that, he can maximize his profit by minimizing customer blocking. The strategy of the \emph{leader} is to tune the pricing parameter
vector $\boldsymbol{\theta}$ at each station. Then the offered $1\times N$ price vector becomes a function of vector $\theta$, $\boldsymbol{p}_k=h(\boldsymbol{\theta})$.
Also we assume that regulations and policies do forbids the network operator to take advantage of the EV's location information to ask for a high price of those
whose battery levels are very low. Similarly network operator can discard customers who tries to cheat by re-entering a congested station many times. The utility function of the network operator becomes:
\begin{align}\label{leader_util}
& {\mathop {\arg \max}\limits_{\boldsymbol{\theta}}}\;  &&\boldsymbol{p}_k(\theta)\boldsymbol{q}-p_{{\!_B}}\boldsymbol{q}_{\!_B}&\\
&\text{s.t.}&&\boldsymbol{p}_k,\boldsymbol{q},\;\boldsymbol{q}_{\!_B}\in\mathbb{R^N},\;p_{\!_B}\in\mathbb{R^{+}}
\end{align}
where, $\boldsymbol{q}$ is calculated from~\eqref{customer_util} and equals to $\boldsymbol{e}_n$. Similarly,
\begin{align}
\boldsymbol{q_{\!_B}} = \left\{ {\begin{array}{*{20}{c l}}
\boldsymbol{0},&\mbox{if gets service}\\
\boldsymbol{e}_n,&\mbox{if blocked.}
\end{array}} \right.
\end{align}
and the price vector $\boldsymbol{p}_k$ is the same as in \eqref{customer_util}. Finally, $p_{\!_B}$ is a scalar and represents priced paid when
blocking occurs. The algorithm for solving the decentralized control is presented in Algorithm~\ref{algo1}.
 \subsection{Applicability of the Model}
In addition to load balancing in a network of fast charging stations, the proposed model is relevant to a variety of applications. In a broader context, the model can be applied to minimize the cost of wide range of scenarios in contestable network composed of selfish users with hierarchy in their decisions (e.g., assignment problem in content distribution networks and job scheduling in cloud computing). In the case of smart grid, the supply and demand scheduling problem for unsplittable loads can be formulated by a hierarchical game, where the owner of the generators act as the leader and the consumers act as the followers. More specifically, the followings can be considered as possible applications.

Consider a demand response management scenario where multiple electric utility companies and consumers interact each other, and the goal of the each entity is to maximize its own benefit. Assume that the cost of the utilities is mainly determined by the power generation method (e.g., unit commitment vs. economic dispatch), associated cost for GHG emissions, and the risk of power outages due to congestions, especially in transmission network. In this case, a network operator can set prices to influence the customers served by high operating cost utilities to other utility options. This will allow consumers to receive discounted services for their deferrable and unsplittable demand (e.g., washing machine). A similar scenario can be applied to energy trading among microgrids and the minimum cost assignments can be orchestrated via \emph{leader (sellers)-follower (buyers)} game.

\section{Numerical Results}

\subsection{Single Charging Station}
\begin{algorithm}[t]
\begin{small}
\caption{Decentralized Control}\label{algo1}
\begin{algorithmic}
\REQUIRE $\boldsymbol{\theta} \geq 0,\;
\left|\mathcal{N}\right|,\left|\mathcal{K}\right|\in
\mathbb{Z^{+}}$ \FOR {customer-k $\leftarrow$ $1$ to $K$}  \STATE network owner offers
$\boldsymbol{p}_k(\boldsymbol\theta)$ $\in$ $\mathbb{R^N}$ to
Eq.(\ref{customer_util}) \STATE calculate utility $\boldsymbol{U}_k$
$\in$ $\mathbb{R^N}$ \STATE pick station
$n$=indexOf(min($\boldsymbol{U}_k$)) \STATE enter station, set
$\boldsymbol{e}(n)$=$1$ \IF {gets service} \STATE set
$\boldsymbol{q}$=$\boldsymbol{e}(n)$, $\boldsymbol{q_B}$=$ \boldsymbol{0}$ \ELSE \STATE set
$\boldsymbol{q}$=$\boldsymbol{0}$, $\boldsymbol{q_B}$=$\boldsymbol{e}(n)$ \ENDIF \STATE
calculate Eq.(\ref{leader_util})
 \ENDFOR
\end{algorithmic}
\end{small}
\end{algorithm}

The motivation for the proposed control scheme is based on the dynamics of the underlying charging station model presented in section \ref{preWork}.
In this section, we present numerical examples to highlight the relationship between the charging station serving capacity (grid power $S$ and storage unit $R$), customer demand ($\lambda(t)$), and the station policy (QoS). Technological constraints for the EV charging rate and the local storage unit are held constant throughout the simulations and set as $\mu$=$2$ and $\nu$=$4$. A detailed discussion on the performance of different charging technologies and computation techniques are presented in \cite{sgc12}. In Fig.~\ref{singleBlocking}, the blocking probability ($P_{BT}$) is evaluated for different set of station parameters (by solving the steady state probability distribution of the underlying Markov chain, see \cite{jsac}). Note that as the station resources increase, more customers can be served. However, the power grid limits the serving resources drawn from the grid and the customer demand ($\lambda$) is mainly determined by the location of the charging station and the time of the day. Therefore, the idea behind the proposed customer routing scheme is to off-load the traffic to neighboring stations when the grid resources are inadequate. This is further clarified in the next two examples.

 In the first one, we present the profit model (details are given in section \ref{preWork}) for a single charging station for varying customer demand. We assume that each customer pays $p_{normal}$=$4$ units of money for the charging service. On the other hand, if a customer is blocked, station operator pays a penalty $p_B$=$5$ units of money. As explained before the penalty motivates station operator to provide a good QoS level. The results presented in Fig.~\ref{singleProfit} shows that there is an optimal traffic regime for individual stations for a given set of station parameters. Note that under light traffic ($\lambda$=[$1$,$2$] )station is making a negative profit mainly because the energy storage acquisition cost and the cost of grid power outweighs the revenue gained from the customers. Similarly, under heavy traffic the cost of customer blockings decreases the station revenue. To that end, the proposed scheme aims to balance the load among the station network by routing customers to idle neighboring stations.

 Another goal of the routing scheme is to ensure QoS levels at each charging stations. The QoS targets are mainly determined by the station operator, and whenever the customer demand exceeds the blocking target, station operator employs pricing-based routing scheme to balance the load. Suppose that station has the following resources $S$=$5$ and $R$=$5$ and for sake of simplicity arrival rate follows a sine wave $\lambda$=$5$+$1.75\sin (2\pi t/80)$. Furthermore, the station policy for the QoS is set to $1\%$, meaning that station operator wants to serve $99\%$ of customers. Then, as depicted in Fig. \ref{sineWave}, the proposed routing scheme is applied during the congestion periods from time slot $22$ to $58$.

\subsection{Network of Charging Stations}
In this part, we illustrate the proposed resource allocation and EV routing scheme.
The simulation scenario is set as follows: the charging network comprises of five stations in a $30\times 30$ unit square area. The coordinates of the locations of them are: $(5, 25)$, $(10, 10)$, $(25, 25)$, $(15, 15)$, and $(25, 5)$ for stations 1 to 5 respectively. To capture the spatial variability of EVs we used the following mixture distribution, whose parameters are calibrated based on results from \cite{jsac}.
\begin{equation}\label{spatDistX}
f(X) =
\begin{cases}
    15\times Be(4.42, 0.763), & \text{0 $\le$ X $\le$ 15} \\
\end{cases}
\end{equation}
\begin{equation}\label{spatDistY}
f(Y) =
\begin{cases}
    15\times Be(2.42, 0.799), & \text{0 $\le$ Y $\le$ 15} \\
\end{cases}
\end{equation}
where $Be$ denotes the Beta distribution function. The spatial
distribution and the locations of the stations are given in
Fig.~\ref{mapLayout}. As depicted in this figure, half of the
EV population resides in the lower left of the area (e.g.,
downtown region) under consideration, while the remaining
half in the rest of the area. Given this spatial distribution of
customers and in the absence of a decentralized control
mechanism, EV demand for each charging node would be $1\%, 50\%,
42\%, 2\%$, and $5\%$ for stations $1$-$5$, respectively. As a
baseline scenario we assume that no allocation or customer
routing of any kind takes places. To begin with, let us assume
that the utility company can provide $S^{max}$=$39$ units of power
and the distribution grid constraint $S^{limit}$=$13$. for the
baseline scenario stations get
\begin{small}$\boldsymbol{S}$=$[7$,$\;8$,$\;8$,$\;8$,$\;8]$\end{small}
respectively. In this case the blocking probability performance
of five stations would be
\begin{small}$\boldsymbol{B}$=$[ \mathtt{\sim}0,\;0.46,\;0.36,\;\mathtt{\sim}0,\; \mathtt{\sim}0]$\end{small}.
This leads to a highly undesirable situation; stations $2$ and
$3$ performs a very bad service, whereas neighboring stations
are over-provisioned and the serving resources are wasted. Our
two step framework, optimal resource allocation and customer
routing will improve the system performance, hence the
percentage of served customers significantly. Details are
presented next.

Given the customer demand for each station as the baseline scenario, the serving capacity is allocated among the stations according to equations \eqref{AllocateS} and \eqref{excessive}. Then the optimal resource allocation would be \begin{small}$\boldsymbol{B}$=$[6,\;13,\;13,\;3,\;4]$\end{small}. Note that due to power grid constraints, the excess power $S_2^{excess}$=$5$ and $S_3^{excess}$=$2$ is allocated to stations 1, 4, and 5 as $[5$,$1$,$1]$. Moreover, the charge rate to satisfy one EV charging request is $\mu_1$= $\dots$ =$\mu_5$=$2$, while the charging rate from local energy storage unit is $\nu_1$= $\dots$ =$\nu_5$=$3$. Also, each station employs an energy storage size of $R_i=8$, $\forall i\in\mathcal{N}$.

The details of the discrete event simulation is presented next. The network-wide charging requests (or arrival rate) is depicted in Fig.~\ref{custDemand}. EV demand at each station is proportional to the traffic load calculated above; for instance, for station $2$ it is
about $\lambda_2$($t$=$10am$)=$52\times 0.5$ = $26$ (EVs arrive in one hour). The charging station operator aims to provide service to the EVs with QoS guarantee
$\delta$ = $0.05$ (at least $95\%$ of the customers will be served at all times) and dissatisfaction parameter $\xi$=$0.1$ for all customers. Based on these charging station parameters, the \emph{leader} (system operator)
solves \eqref{optLambda} and initiates the routing game only when the arrival rate exceeds this threshold, which corresponds to the
following thresholds for each individual station: $\boldsymbol{\lambda_i^{*}}=[6.7,\;23.4,\;23.4,\;3.3,\;5.0] $. Since it may be challenging (and possibly wasteful) to update arrival rates in real time, we set $15$ minute intervals at which arrival rates are recalculated. Just to clarify how resource allocation improves the system performance, let us compare it with baseline scenario at 4pm according to weighted sum of overall system performance using the formula given below:
\begin{equation}
\sum\limits_{i \in l} {{w_i}{B_i}},\; \text{where}\;{w_i} = \frac{{{\lambda _i}}}{{\sum\limits_{i \in l} {{\lambda _i}} }}
\end{equation}
where $\lambda_i$=$\tilde{\lambda}_{\!_{EV}} $ Then, the weighted sum of the blocking probabilities for the baseline scenario and the power resource allocation would be $0.38$ and $0.10$ which leads to more efficient use of grid resources.

The details of the customer routing framework is discussed next. In the absence of congestion
($\lambda_i< \lambda_i^*$), customers pay the operator $p_{normal}$=$4$ units,
whereas if EVs are blocked the operator rewards them with $p_{block}$=$6$. As
mentioned before, this cost is a penalty to the operator for poor service
which can impact customer loyalty and its long term reputation. Also, to
calculate blocking probabilities we set $\gamma _1$=$0.45$ and $\gamma
_2$=$0.55$. For the $k^{th}$ EV, $c_{k_{inctv}}$ is a uniformly distributed
random number in the interval $[0.75, 1.0]$ and $p_{dis}$ is a uniformly
distributed random number (per unit of distance) in the interval $[0.02,
0.05]$. We assume that driving duration is linearly correlated with distance
(based on an average speed of $40$ mph). Currently popular EV models (e.g.,
Nissan Leaf) exhibit $0.22$ kWh/mile energy consumption; thus, we set
$p_{drive}$ = $0.03$ per unit of distance.

Next, the strategy of the \emph{leader} (setting the tuning
parameter vector $\boldsymbol{\theta}$) is given. As stations $1$, $4$, and
$5$ operate under light traffic, the operator is motivated shape
excessive demand at stations $2$ and $3$ and route customers to
these stations. For the given arrival rate at
($\lambda(t=4pm)$), we set
$\theta_1$=$\theta_4$=$\theta_5$=$0.5$ and investigated how
$\theta_2$ and $\theta_3$ affects the blocking probability at
these
 stations. As depicted in Fig.~\ref{leaderStrategy} station $2$ and $3$ violate $\delta_{max}$=$0.05$ QoS policy, however station operator increases $\theta_2$ and $\theta_3$ to set the congestion prices and start routing customers to neighboring stations.
 Our customer routing framework enables stations to meet SLA agreements. Note that $\theta_2$ and $\theta_3$>$0.4$ does not lead to further reduction in blocking probability because the neighboring stations also hit $\lambda^*$ and does not accept more customers. Similarly, Fig.~\ref{loadBalancing} shows the load balancing framework. Notice that for low customer demand there is no need for customer routing scheme because the regional serving capacity resources suffice the demand. On the other hand as stations become more congested customers are routed to neighboring stations. One important point of consideration is that since station $1$ is physically more close to congested stations $2$ and $3$, more customers accept (almost $11$\% of the population) to go to this station. In Fig. \ref{IncreasedServings} we compute the percentage of increase in the number of served customers respect to base line scenario. For the given simulation setting, the proposed framework serves close $50\%$ more customers with the same amount of grid resources.

The proposed decentralized control mechanism is simulated with
the given set of parameters and compared to baseline and power
allocation scenarios. The performance of the most congested
station (station $2$) is presented in Fig.~\ref{station2}. It
can be seen that the percentage of served customers can be
increased significantly by employing proposed framework. We
further compare the network-wide net revenue. As depicted
Fig.~\ref{netRev}, the network operator optimally allocates grid
resources and do load balancing, the utilization
 of each charging station will increase, and hence more customers will be served with the same amount of grid resources. Hence, by employing proposed framework network operator increases her profit tremendously.
\vspace{-10pt}
\section{Conclusion}
In this paper, we introduce a Stackelberg game theoretic based control mechanism to manage a population of self-interested mobile EVs.
The aim of the \emph{leader} (Network Operator) is to serve more customers with the same amount of grid resources, while the goal of the EV \emph{followers} (EV drivers) is to get the charging service at a minimum cost. Utility functions are developed to represent the behavior of both parties.

 As a future work, we are aiming to investigate the communication components involved in making this scheme operational. Charging stations communicate with customers via smart applications using 3G/4G wireless technologies. The successful deployment of the coordination mechanism discussed in this paper, heavily depends on the availability of the necessary communication infrastructure to ensure reliable information dissemination. In our framework, mobile EVs use wireless communication technologies to locate and
retrieve pricing information. On the other hand, wireless communications enable the network operator to monitor station
usage, update prices and forward them to EVs. Hence, it is crucial to identify the communication system requirements (e.g., communication losses, delays etc.) and
quantify its impact on the performance of the charging stations.

\bibliographystyle{IEEEtran}
\bibliography{tsg}
\vspace{-35pt}
 \begin{IEEEbiographynophoto}{I. Safak Bayram}
(S10, M14)
received the B.S. degree in electrical and electronics engineering from Dokuz Eylul University, Izmir, Turkey in 2007, the M.S. degree in Telecommunications from the University of Pittsburgh
in 2010, and the Ph.D. degree in Computer Engineering from North Carolina State University, in 2014. He received the Best Paper Award at the Third IEEE International Conference on Smart Grid Communications and the Student Travel Grant at a previous Smart Grid Communications Conference.
Currently he is a Postdoctoral Research Scientist in the Department of Electrical and Computer Engineering at Texas A\&M University at Qatar. His research interest include optimization, control, and stochastic modeling of communications and power networks.
 \end{IEEEbiographynophoto}
\vspace{-40pt}
\begin{IEEEbiographynophoto}{George Michailidis}
received the Ph.D. degree in
mathematics from the University of California, Los
Angeles, in 1996.
He was a Postdoctoral Fellow in the Department
of Operations Research at Stanford University from
1996 to 1998. He joined the University of Michigan,
Ann Arbor, in 1998, where he is currently a Professor
of Statistics, Electrical Engineering, and Computer
Science. He is a Fellow of the Institute of Mathematical Statistics, the American
Statistical Association and the International Statistical Institute. He is the editor-in-chief
of the Electronic Journal of Statistics and serves on a number of editorial boards. He served as
symposium chair on
Support for Storage, Renewable Sources and Micro-grid for SmartGridComm 2012 and is
the secretary of the IEEE subcommittee on SmartGrid Communications.
His research interests are in the areas of stochastic
network modeling and performance evaluation,
queuing analysis and congestion control, statistical modeling
and analysis of Internet traffic, network tomography, and analysis
of high dimensional data with network structure.
 \end{IEEEbiographynophoto}
\vspace{-40pt}
\begin{IEEEbiographynophoto}{Michael Devetsikiotis}
(S85, M94, SM03, F12)  received the
Dipl.Ing. degree in electrical engineering from the
Aristotle University of Thessaloniki, Greece, in
1988, and the M.Sc. and Ph.D. degrees in electrical
engineering from North Carolina State University,
Raleigh, in 1990 and 1993, respectively.
In 1993 he joined the Broadband Networks Laboratory
at Carleton University, Ottawa, ON, Canada,
as a Postdoctoral Fellow. He later became an Adjunct
Research Professor in the Department of Systems and
Computer Engineering at Carleton University in 1995, an Assistant Professor in
1996 and an Associate Professor in 1999. He joined the Department of Electrical
and Computer Engineering at North Carolina State as an Associate Professor,
in 2000, and became a Professor in 2006. He served as Chairman of the IEEE Communications
Society Technical Committee Communication Systems Integration
and Modeling, and as a member of the ComSoc Education Board. He has
also served as an Associate Editor of the IEEE Communications Letters, an Area
Editor of the ACM Transactions on Modeling and Computer Simulation, and
on the editorial boards of the International Journal of Simulation and Process
Modeling, the IEEE Communications Surveys and Tutorials, and the Journal
of Internet Engineering.
 \end{IEEEbiographynophoto}

\end{document}